\documentclass[11pt]{article}

\usepackage{fullpage}
\usepackage{distopt}

\renewcommand{\graph}{\ensuremath{\mc{G}}}

\def\extrastep{\ensuremath{\eta}}

\def\tcycle{\ensuremath{T_{\rm cycle}}}
\def\tcent{\ensuremath{T_{\rm cent}}}
\def\tdist{\ensuremath{T_{\rm dist}}}
\def\mdist{\ensuremath{m_{\rm dist}}}
\def\mcycle{\ensuremath{m_{\rm cycle}}}
\renewcommand{\order}{\mathcal{O}}

\def\avgx{\ensuremath{\what{x}}}
\def\xavg{\avgx}

\def\diam{D}

\def\comcost{\mathsf{C}}
\def\comcost{\ensuremath{\theta}}
\def\comcost{\ensuremath{C}}

\makeatletter
\long\def\@makecaption#1#2{
  \vskip 0.8ex
  \setbox\@tempboxa\hbox{\small {\bf #1:} #2}
  \parindent 1.5em  
  \dimen0=\hsize
  \advance\dimen0 by -3em
  \ifdim \wd\@tempboxa >\dimen0
  \hbox to \hsize{
    \parindent 0em
    \hfil 
    \parbox{\dimen0}{\def\baselinestretch{0.96}\small
      {\bf #1.} #2
    } 
    \hfil}
  \else \hbox to \hsize{\hfil \box\@tempboxa \hfil}
  \fi
}
\makeatother

\title{Distributed Stochastic Optimization with Centralized Control}
\author{Alekh Agarwal \and John C. Duchi \\
{\tt \{alekh,jduchi\}@eecs.berkeley.edu}}

\def\comment#1{}

\begin{document}

\begin{center}
  {\bf {\Large Distributed Delayed Stochastic Optimization}} \\
  \vspace{.25cm}
  Alekh Agarwal ~~~~~~~~~~ John C. Duchi \\
  \vspace{.1cm}
  {\tt \{alekh,jduchi\}@eecs.berkeley.edu} \\
  \vspace{.2cm}
  Department of Electrical Engineering and Computer Sciences \\
  University of California, Berkeley
\end{center}

\begin{abstract}
  We analyze the convergence of gradient-based optimization algorithms
  that base their updates on delayed stochastic gradient
  information. The main application of our results is to the
  development of gradient-based distributed optimization algorithms
  where a master node performs parameter updates while worker nodes
  compute stochastic gradients based on local information in parallel,
  which may give rise to delays due to asynchrony. We take motivation
  from statistical problems where the size of the data is so large
  that it cannot fit on one computer; with the advent of huge datasets
  in biology, astronomy, and the internet, such problems are now
  common. Our main contribution is to show that for smooth stochastic
  problems, the delays are asymptotically negligible and we can
  achieve order-optimal convergence results. In application to
  distributed optimization, we develop procedures that overcome
  communication bottlenecks and synchronization requirements. We show
  $n$-node architectures whose optimization error in stochastic
  problems---in spite of asynchronous delays---scales asymptotically
  as $\order(1 / \sqrt{nT})$ after $T$ iterations. This rate is known
  to be optimal for a distributed system with $n$ nodes even in the
  absence of delays. We additionally complement our theoretical
  results with numerical experiments on a statistical machine learning
  task.
\end{abstract}

\section{Introduction}

We focus on stochastic convex optimization problems of the form
\begin{equation}
  \minimize_{x \in \xdomain} ~ f(x)
  ~~~ {\rm for} ~~~
  f(x) \defeq \E_P[F(x; \statsample)]
  = \int_\statdomain F(x; \statsample) dP(\statsample),
  \label{eqn:stochastic-objective}
\end{equation}
where $\xdomain \subseteq \R^d$ is a closed convex set, $P$ is a
probability distribution over $\statdomain$, and $F(\cdot\,;
\statsample)$ is convex for all $\statsample \in \statdomain$, so that
$f$ is convex. The goal is to find a parameter $x$ that approximately
minimizes $f$ over $x \in \xdomain$.  Classical stochastic gradient
algorithms~\cite{RobbinsMo51,Polyak87} iteratively update a parameter
$x(t) \in \xdomain$ by sampling $\statsample \sim P$, computing $g(t)
= \nabla F(x(t); \statsample)$, and performing the update $x(t + 1) =
\Pi_{\xdomain}(x(t) - \stepsize(t) g(t))$, where $\Pi_{\xdomain}$
denotes projection onto the set $\xdomain$.  In this paper, we analyze
asynchronous gradient methods, where instead of receiving current
information $g(t)$, the procedure receives out of date gradients $g(t
- \tau(t)) = \nabla F(x(t - \tau(t)), \statsample)$, where $\tau(t)$
is the (potentially random) delay at time $t$.  The central
contribution of this paper is to develop algorithms that---under
natural assumptions about the functions $F$ in the
objective~\eqref{eqn:stochastic-objective}---achieve asymptotically
optimal rates for stochastic convex optimization in spite of delays.

Our model of delayed gradient information is particularly relevant in
distributed optimization scenarios, where a master maintains the
parameters $x$ while workers compute stochastic gradients of the
objective~\eqref{eqn:stochastic-objective}. The architectural
assumption of a master with several worker nodes is natural for
distributed computation, and other researchers have considered models
similar to those in this paper~\cite{NedicBeBo01,LangfordSmZi09}. By
allowing delayed and asynchronous updates, we can avoid
synchronization issues that commonly handicap distributed systems.

Certainly distributed optimization has been studied for several
decades, tracing back at least to seminal work of Tsitsiklis and
colleagues~(\cite{Tsitsiklis84,BertsekasTs89}) on minimization of
smooth functions where the parameter vector is distributed. More
recent work has studied problems in which each processor or node $i$
in a network has a local function $f_i$, and the goal is to minimize
the sum $f(x) = \ninv \sum_{i=1}^n
f_i(x)$~\cite{NedicOz09,RamNeVe10,JohanssonRaJo09,DuchiAgWa10}. Most
prior work assumes as a \emph{constraint} that data lies on several
different nodes throughout a network.  However, as Dekel et
al.~\cite{DekelGiShXi10} first noted, in distributed stochastic
settings independent realizations of a stochastic gradient can be
computed concurrently, and it is thus possible to obtain an aggregated
gradient estimate with lower variance. Using modern stochastic
optimization algorithms (e.g.~\cite{JuditskyNeTa08,Lan10}), Dekel et
al.\ give a series of reductions to show that in an $n$-node network
it is possible to achieve a speedup of $\order(n)$ over a
single-processor so long as the objective $f$ is smooth.

\begin{figure}[t]
  \begin{center}
    \includegraphics[width=.65\columnwidth]{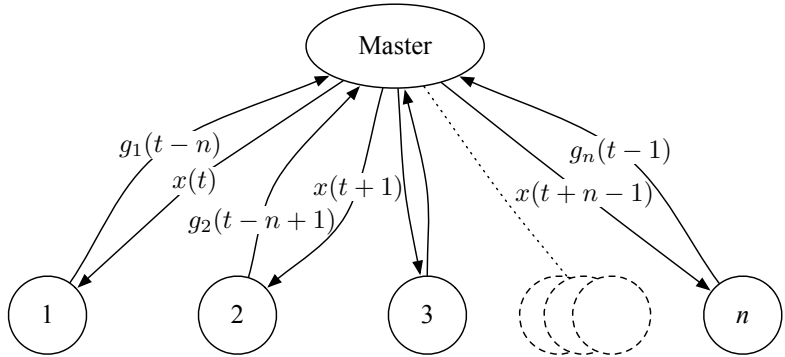}
    \caption{\label{fig:master-worker} Cyclic delayed update
      architecture. Workers compute gradients cyclically and in parallel,
      passing out-of-date information to master. Master responds with current
      parameters. Diagram shows parameters and gradients communicated
      between rounds $t$ and $t+n-1$.}
  \end{center}
\end{figure}

Our work is closest to Nedi\'c et al.'s asynchronous subgradient
method~\cite{NedicBeBo01}, which is an incremental gradient procedure
in which gradient projection steps are taken using out-of-date
gradients. See Figure~\ref{fig:master-worker} for an illustration. The
asynchronous subgradient method performs non-smooth minimization and
suffers an asymptotic penalty in convergence rate due to the delays:
if the gradients are computed with a delay of $\tau$, then the
convergence rate of the procedure is $\order(\sqrt{\tau / T})$.
The setting of distributed optimization provides an elegant
illustration of the role played by the delay in convergence rates. As
in Fig.~\ref{fig:master-worker}, the delay $\tau$ can essentially be
of order $n$ in Nedi\'c et al.'s setting, which gives a convergence
rate of $\order(\sqrt{n/T})$. A simple centralized stochastic gradient
algorithm attains a rate of $\order(1/\sqrt{T})$, which suggests
something is amiss in the distributed algorithm.
Langford et al.~\cite{LangfordSmZi09} rediscovered Nedi\'c et al.'s
results and attempted to remove the asymptotic penalty by considering
smooth objective functions, though their approach has a technical
error (see Appendix~\ref{appendix:langford-error}), and even so they
do not demonstrate any provable benefits of distributed
computation. We analyze similar asynchronous algorithms, but we show
that for smooth stochastic problems the delay is asymptotically
negligible---the time $\tau$ does not matter---and in fact, with
parallelization, delayed updates can give provable performance
benefits.



We build on results of Dekel et al.~\cite{DekelGiShXi10}, who show
that when the objective $f$ has Lipschitz-continuous gradients, then
when $n$ processors compute stochastic gradients in parallel using a
common parameter $x$ it is possible to achieve convergence rate
$\order(1 / \sqrt{Tn})$ so long as the processors are synchronized
(under appropriate synchrony conditions, this holds nearly
independently of network topology). A variant of their approach is
asymptotically robust to asynchrony so long as most processors remain
synchronized for most of the time~\cite{DekelGiShXi10b}. We show
results similar to their initial discovery, but we analyze the effects
of asynchronous gradient updates where all the nodes in the network
can suffer delays. Application of our main results to the distributed
setting provides convergence rates in terms of the number of nodes $n$
in the network and the stochastic process governing the
delays. Concretely, we show that under different assumptions on the
network and delay process, we achieve convergence rates ranging from
$\order(n^3/T + 1/\sqrt{Tn})$ to $\order(n/T + 1 / \sqrt{Tn})$, which
is $\order(1/\sqrt{nT})$ asymptotically in $T$. For problems with large $n$,
we demonstrate faster rates ranging from
$\order((n/T)^{2/3} + 1/\sqrt{Tn})$ to $\order(1/T^{2/3} + 1 /
\sqrt{Tn})$. In either case, the time necessary to achieve
$\epsilon$-optimal solution to the
problem~\eqref{eqn:stochastic-objective} is asymptotically $\order(1 /
n\epsilon^2)$, a factor of $n$---the size of the network---better than
a centralized procedure in spite of delay.

The remainder of the paper is organized as follows. We begin by reviewing
known algorithms for solving the stochastic optimization
problem~\eqref{eqn:stochastic-objective} and stating our main assumptions.
Then in Section~\ref{sec:main-results} we give abstract descriptions of our
algorithms and state our main theoretical results, which we make concrete in
Section~\ref{sec:distributed} by formally placing the analysis in the setting
of distributed stochastic optimization. We complement the theory in
Section~\ref{sec:experiments} with experiments on a real-world dataset, and
proofs follow in the remaining sections.

\paragraph{Notation}
For the reader's convenience, we collect our (mostly standard) notation
here. We denote general norms by $\norm{\cdot}$, and the dual norm
$\dnorm{\cdot}$ to the norm $\norm{\cdot}$ is defined as $\dnorm{z} \defeq
\sup_{x : \norm{x} \le 1} \<z, x\>$.
The subdifferential set of a function $f$ is
\begin{equation*}
  \partial f(x)
  \defeq \left\{g \in \R^d \mid
  f(y) \ge f(x) + \<g, y - x\> ~ \mbox{for all }y \in \dom f\right\}
\end{equation*}
We use the shorthand $\dnorm{\partial f(x)} \defeq \sup_{g \in
  \partial f(x)} \dnorm{g}$. A function $f$ is $G$-Lipschitz with
respect to the norm $\norm{\cdot}$ on $\xdomain$ if for all $x, y \in
\xdomain$, $|f(x) - f(y)| \le G\norm{x - y}$. For convex $f$, this is
equivalent to $\dnorm{\partial f(x)} \le G$ for all $x \in \xdomain$
(e.g.~\cite{HiriartUrrutyLe96}). A function $f$ is $L$-smooth on
$\xdomain$ if $\nabla f$ is Lipschitz continuous with respect to the
norm $\norm{\cdot}$, defined as
\begin{equation*}
  \dnorm{\nabla f(x) - \nabla f(y)} \le L \norm{x - y},
  ~~~ \mbox{equivalently}, ~~~
  f(y) \le f(x) + \<\nabla f(x), y - x\> + \frac{L}{2} \norm{x - y}^2.
\end{equation*}
For convex differentiable $h$, the Bregman
divergence~\cite{Bregman67} between $x$ and $y$ is defined as
\begin{equation}
  \label{eqn:def-bregman-divergence}
  D_h(x, y) \defeq h(x) - h(y) - \<\nabla h(y), x - y\>.
\end{equation}
A convex function $h$ is $c$-strongly convex with
respect to a norm $\norm{\cdot}$ over $\xdomain$ if
\begin{equation}
  \label{eqn:strong-convex}
  h(y) \ge h(x) + \<g, y - x\>
  + \frac{c}{2} \norm{x - y}^2
  ~~~ \mbox{for~all}~
  x, y \in \xdomain ~ \mbox{and} ~ g \in \partial h(x).
\end{equation}
We use $[n]$ to denote the set of integers $\{1, \ldots, n\}$.

\section{Setup and Algorithms}
\label{sec:setup}

In this section we set up and recall the delay-free algorithms
underlying our approach. We then give the appropriate delayed
versions of these algorithms, which we analyze in the sequel.

\subsection{Setup and Delay-free Algorithms}
\label{sec:setup-cent}

To build intuition for the algorithms we analyze, we first describe
two closely related first-order algorithms: the dual averaging
algorithm of Nesterov~\cite{Nesterov09} and the mirror descent
algorithm of Nemirovski and Yudin~\cite{NemirovskiYu83}, which is
analyzed further by Beck and Teboulle~\cite{BeckTe03}. We begin by
collecting notation and giving useful definitions. Both algorithms are
based on a proximal function $\prox(x)$, where it is no loss of
generality to assume that $\prox(x) \ge 0$ for all $x \in \xdomain$.
We assume $\prox$ is $1$-strongly convex (by scaling, this is no loss
of generality). By
definitions~\eqref{eqn:def-bregman-divergence}
and~\eqref{eqn:strong-convex}, the divergence $\divergence$ satisfies
$\divergence(x, y) \ge \half \norm{x - y}^2$.

In the oracle model of stochastic optimization that we assume, at time
$t$ both algorithms query an oracle at the point $x(t)$, and the
oracle then samples $\statsample(t)$ i.i.d.\ from the distribution $P$
and returns $g(t) \in \partial F(x(t); \statsample(t))$. The dual
averaging algorithm~\cite{Nesterov09} updates a dual vector $z(t)$ and
primal vector $x(t) \in \xdomain$ via
\begin{equation}
  \label{eqn:da-update}
  z(t + 1) = z(t) + g(t)
  ~~~ \mbox{and} ~~~
  x(t + 1) = \argmin_{x \in \xdomain} \Big\{
  \<z(t + 1), x\> + \frac{1}{\stepsize(t + 1)}\prox(x)\Big\},
\end{equation}
while mirror descent~\cite{NemirovskiYu83,BeckTe03} performs the update
\begin{equation}
  \label{eqn:md-update}
  x(t + 1) = \argmin_{x \in \xdomain} \Big\{
  \<g(t), x\> + \frac{1}{\stepsize(t)} \divergence(x, x(t))\Big\}.
\end{equation}
Both make a linear approximation to the function being minimized---a
global approximation in the case of the dual averaging
update~\eqref{eqn:da-update} and a more local approximation for mirror
descent~\eqref{eqn:md-update}---while using the proximal function
$\prox$ to regularize the points $x(t)$.

We now state the two essentially standard
assumptions~\cite{JuditskyNeTa08,Lan10,Xiao10} we most often make
about the stochastic optimization
problem~\eqref{eqn:stochastic-objective}, after which we recall the
convergence rates of the algorithms~\eqref{eqn:da-update}
and~\eqref{eqn:md-update}.
\begin{assumption}[Lipschitz Functions]
  \label{assumption:lipschitz-func}
  For $P$-a.e.\ $\statsample$, the function $F(\cdot\,; \statsample)$ is
  convex. Moreover, for any $x \in \xdomain$, $\E[\dnorm{\partial F(x;
      \statsample)}^2] \le G^2$.
\end{assumption}
\noindent
In particular, Assumption~\ref{assumption:lipschitz-func} implies that
$f$ is $G$-Lipschitz continuous with respect to the norm
$\norm{\cdot}$ and that $f$ is convex. Our second assumption has been
used to show rates of convergence based on the variance of a gradient
estimator for stochastic optimization problems
(e.g.~\cite{JuditskyNeTa08, Lan10}).
\begin{assumption}[Smooth Functions]
  \label{assumption:lipschitz-grad}
  The function $f$ defined in~\eqref{eqn:stochastic-objective} has
  $L$-Lipschitz continuous gradient, and for all $x \in \xdomain$ the
  variance bound $\E[\dnorm{\nabla f(x) - \nabla F(x; \statsample)}^2]
  \le \stddev^2$ holds.\footnote{If $f$ is differentiable, then
    $F(\cdot; \statsample)$ is differentiable for
    $P$-a.e.\ $\statsample$, and conversely, but $F$ need not be
    smoothly differentiable~\cite{Bertsekas73}. Since $\nabla F(x;
    \statsample)$ exists for $P$-a.e.\ $\statsample$, we will write
    $\nabla F(x; \statsample)$ with no loss of generality.}
\end{assumption}
\noindent
Several commonly used functions satisfy the above assumptions, for example:
\begin{enumerate}[(i)]
\item The \emph{logistic loss}: $F(x; \statsample) = \log[1 +
  \exp(\<x,\statsample\>)]$, the objective for logistic regression in
  statistics (e.g.~\cite{HastieTiFr01}). The objective $F$ satisfies
  Assumptions~\ref{assumption:lipschitz-func}
  and~\ref{assumption:lipschitz-grad} so long as $\norm{\statsample}$
  is bounded.
  \label{item:logistic-regression}
\item \emph{Least squares} or \emph{linear regression}: $F(x;\statsample) = (a
  - \<x,b\>)^2$ where $\statsample = (a,b)$ for $a \in \R^d$ and $b \in \R$,
  satisfies
  Assumptions~\ref{assumption:lipschitz-func}
  and~\ref{assumption:lipschitz-grad} as long as $\statsample$ is bounded and
  $\xdomain$ is compact.
  \label{item:least-squares}
\end{enumerate}
\noindent
We also make a standard compactness assumption on the optimization set
$\xdomain$. 
\begin{assumption}[Compactness]
  For $x^* \in \argmin_{x \in \xdomain} f(x)$ and $x \in \xdomain$, 
  the bounds $\prox(x^*)
  \le \radius^2/2$ and \mbox{$\divergence(x^*,x) \leq \radius^2$} both hold.
  \label{assumption:compact}
\end{assumption}

Under Assumptions~\ref{assumption:lipschitz-func}
or~\ref{assumption:lipschitz-grad} in addition to
Assumption~\ref{assumption:compact}, the updates~\eqref{eqn:da-update}
and~\eqref{eqn:md-update} have known convergence rates. Define the
time averaged vector $\avgx(T)$ as
\begin{equation}
  \label{eqn:avg-x}
  \avgx(T) \defeq \frac{1}{T} \sum_{t=1}^T x(t + 1).
\end{equation}
Then under Assumption~\ref{assumption:lipschitz-func}, both algorithms 
satisfy 
\begin{equation}
  \label{eqn:standard-convergence}
  \E[f(\what{x}(T))] - f(x^*) = \order\left(\frac{\radius
    G}{\sqrt{T}}\right)
\end{equation}
for the stepsize choice $\stepsize(t) = \radius / (G\sqrt{t})$
(e.g.~\cite{Nesterov09,Xiao10,NemirovskiJuLaSh09}).
The result~\eqref{eqn:standard-convergence} is sharp to constant
factors in general~\cite{NemirovskiYu83,AgarwalBaRaWa10}, but can be
further improved under
Assumption~\ref{assumption:lipschitz-grad}. Building on work of
Juditsky et al.~\cite{JuditskyNeTa08} and Lan~\cite{Lan10}, Dekel et
al.~\cite[Appendix A]{DekelGiShXi10} show that under
Assumptions~\ref{assumption:lipschitz-grad}
and~\ref{assumption:compact} the stepsize choice $\stepsize(t)^{-1} =
L + \extrastep(t)$, where $\extrastep(t)$ is a damping factor
set to $\extrastep(t) = \stddev \radius \sqrt{t}$,
yields for either of the updates~\eqref{eqn:da-update} or~\eqref{eqn:md-update}
the convergence rate
\begin{equation}
  \label{eqn:smooth-da-stochastic}
  \E[f(\avgx(T))] - f(x^*) = \order\left(\frac{L\radius^2}{T} +
  \frac{\stddev \radius}{\sqrt{T}}\right).
\end{equation}


\subsection{Delayed Optimization Algorithms}
\label{sec:setupdelayed}

We now turn to extending the dual averaging~\eqref{eqn:da-update} and
mirror descent~\eqref{eqn:md-update} updates to the setting in which instead
descent~\eqref{eqn:md-update} updates to the setting in which instead
of receiving a current gradient $g(t)$ at time $t$, the procedure
receives a gradient $g(t - \tau(t))$, that is, a stochastic gradient
of the objective~\eqref{eqn:stochastic-objective} computed at the
point $x(t - \tau(t))$. In the simplest case, the delays are uniform
and $\tau(t) \equiv \tau$ for all $t$, but in general the delays may
be a non-i.i.d.\ stochastic process. Our analysis admits any sequence
$\tau(t)$ of delays as long as the mapping $t \mapsto \tau(t)$
satisfies $\E[\tau(t)] \le B < \infty$. We also require that each
update happens once, i.e., $t \mapsto t - \tau(t)$ is
one-to-one, though this second assumption is easily satisfied.
                                   
Recall that the problems we consider are stochastic optimization
problems of the form~\eqref{eqn:stochastic-objective}. Under the
assumptions above, we extend the mirror descent and dual averaging
algorithms in the simplest way: we replace $g(t)$ with $g(t - \tau(t))$.
For dual averaging (c.f.\ the update~\eqref{eqn:da-update}) this yields
\begin{equation}
  \label{eqn:da-delayed-update}
  z(t + 1) = z(t) + g(t - \tau(t)) ~~~ \mbox{and} ~~~ x(t + 1) =
  \argmin_{x \in \xdomain} \Big\{\<z(t + 1), x\> +
  \frac{1}{\stepsize(t + 1)} \prox(x) \Big\},
\end{equation}
while for mirror descent (c.f.\ the update~\eqref{eqn:md-update}) we have
\begin{equation}
  \label{eqn:md-delayed-update}
  x(t + 1) = \argmin_{x \in \xdomain} \Big\{ \<g(t - \tau(t)), x\> +
  \frac{1}{\stepsize(t)} \divergence(x, x(t)) \Big\}.
\end{equation}
A generalization of Nedi\'{c} et al.'s results~\cite{NedicBeBo01} by
combining their techniques with the convergence proofs of dual
averaging~\cite{Nesterov09} and mirror descent~\cite{BeckTe03} is
as follows. Under Assumptions~\ref{assumption:lipschitz-func}
and~\ref{assumption:compact}, so long as $\E[\tau(t)] \le B < \infty$ for all
$t$, choosing $\stepsize(t) = \frac{\radius}{G \sqrt{B t}}$ gives rate
\begin{equation}
  \E[f(\what{x}(T))] - f(x^*) = \order\bigg(\frac{\radius G
    \sqrt{B}}{\sqrt{T}}\bigg).
  \label{eqn:nonsmooth-rate}
\end{equation}

\section{Convergence rates for delayed optimization of smooth
  functions} 
\label{sec:main-results}

In this section, we state and discuss several results for asynchronous
stochastic gradient methods. We give two sets of theorems. The first are
for the asynchronous method when we make updates to the parameter vector $x$
using one stochastic subgradient, according to the update
rules~(\ref{eqn:da-delayed-update}) or~(\ref{eqn:md-delayed-update}).
The second method involves using several stochastic subgradients for
every update, each with a potentially different delay, which gives
sharper results that we present in Section~\ref{sec:combo-delays}.

\subsection{Simple delayed optimization}
\label{sec:delayedcyclic}

Intuitively, the $\sqrt{B}$-penalty due to delays for non-smooth
optimization arises from the fact that subgradients can change
drastically when measured at slightly different locations, so a small
delay can introduce significant inaccuracy. To overcome the delay
penalty, we now turn to the smoothness
assumption~\ref{assumption:lipschitz-grad} as well as the Lipschitz
condition~\ref{assumption:lipschitz-func} (we assume both of these
conditions along with Assumption~\ref{assumption:compact} hold for all
the theorems). In the smooth case, delays mean that stale gradients
are only slightly perturbed, since our stochastic algorithms constrain
the variability of the points $x(t)$. As we show in the proofs of the
remaining results, the error from delay essentially becomes a second
order term: the penalty is asymptotically negligible. We study both
update rules~\eqref{eqn:da-delayed-update}
and~\eqref{eqn:md-delayed-update}, and we set $\stepsize(t) =
\frac{1}{L + \extrastep(t)}$. Here $\extrastep(t)$ will be chosen to
both control the effects of delays and for errors from stochastic
gradient information. We prove the following theorem in
Sec.~\ref{sec:stoch-delayed}.
\begin{theorem}
  \label{theorem:da-delayed-stoch}
  Let the sequence $x(t)$ be defined by the update
  \eqref{eqn:da-delayed-update}. Define the stepsize $\extrastep(t)
  \propto \sqrt{t+\tau}$ or let $\extrastep(t) \equiv \extrastep$ for
  all $t$. Then
  \begin{equation*}
    \E \bigg[\sum_{t=1}^T f(x(t + 1))\bigg] - T f(x^*)
    \le \frac{1}{\stepsize(T + 1)} \radius^2
    + \frac{\stddev^2}{2} \sum_{t=1}^T \frac{1}{\extrastep(t)}
    + 2 LG^2(\tau + 1)^2 \sum_{t=1}^T \frac{1}{\extrastep(t-\tau)^2}
    + 2 \tau G\radius.
  \end{equation*}
\end{theorem}
\noindent
The mirror descent update~\eqref{eqn:md-delayed-update} exhibits
similar convergence properties, and we prove the next theorem in
Sec.~\ref{sec:md-stoch-delayed}. 
\begin{theorem}
  \label{theorem:md-delayed-stoch}
  Use the conditions of Theorem~\ref{theorem:da-delayed-stoch} but
  generate $x(t)$ by the update~\eqref{eqn:md-delayed-update}. Then
  \begin{equation*}
    \E\bigg[\sum_{t=1}^T f(x(t + 1))\bigg] - T f(x^*) \le 2 L\radius^2
    + \radius^2 [\extrastep(1) + \extrastep(T)] + \frac{\stddev^2}{2}
    \sum_{t=1}^T \frac{1}{\extrastep(t)} + 2 LG^2 (\tau + 1)^2
    \sum_{t=1}^T \frac{1}{\extrastep(t-\tau)^2} + 2 \tau G \radius
  \end{equation*}
\end{theorem}

In each of the above theorems, we can set $\extrastep(t) =
\stddev\sqrt{t+\tau} / \radius$. As immediate corollaries, we recall the
definition~\eqref{eqn:avg-x} of the averaged sequence of $x(t)$ and use
convexity to see that
\begin{equation*}
  \E[f(\avgx(T))] - f(x^*)
  = \order\left(\frac{L\radius^2 + \tau G\radius}{T}
  + \frac{\stddev \radius}{\sqrt{T}}
  + \frac{LG^2 \tau^2 \radius^2 \log T}{\stddev^2 T}\right)
\end{equation*}
for either update rule.
In addition, we can allow the delay $\tau(t)$ to be
random:
\begin{corollary}
  \label{corollary:random-tau}
  Let the conditions of Theorem~\ref{theorem:da-delayed-stoch}
  or~\ref{theorem:md-delayed-stoch} hold, but allow $\tau(t)$ to be a random
  mapping such that $\E[\tau(t)^2] \le B^2$ for all $t$.  With the choice
  $\extrastep(t) = \stddev \sqrt{T} / \radius$ the
  updates~\eqref{eqn:da-delayed-update} and~\eqref{eqn:md-delayed-update}
  satisfy
  \begin{equation*}
    \E[f(\what{x}(T))] - f(x^*)
    = \order\left(\frac{L\radius^2 + B^2 G\radius}{T}
    + \frac{\sigma \radius}{\sqrt{T}}
    + \frac{LG^2 B^2 \radius^2}{\sigma^2 T}\right).
  \end{equation*}
\end{corollary}
\noindent
We provide the proof of the corollary in
Sec.~\ref{sec:corollary-simple-delayed}. The take-home message from the above
corollaries, as well as Theorems~\ref{theorem:da-delayed-stoch}
and~\ref{theorem:md-delayed-stoch}, is that the penalty in convergence rate
due to the delay $\tau(t)$ is asymptotically negligible. As we discuss
in greater depth in the next section, this has favorable implications
for robust distributed stochastic optimization algorithms.


\subsection{Combinations of delays}
\label{sec:combo-delays}

In some scenarios---including distributed settings similar to those we
discuss in the next section---the procedure has access not to only a single
delayed gradient but to several with different delays.  To abstract away the
essential parts of this situation, we assume that the procedure receives $n$
gradients $g_1, \ldots, g_n$, where each has a potentially different delay
$\tau(i)$.  Now let $\lambda = (\lambda_i)_{i=1}^n$ belong to the probability
simplex, though we leave $\lambda$'s values unspecified for now. Then the
procedure performs the following updates at time $t$: for dual averaging,
\begin{equation}
  z(t+1) = z(t) + \sum_{i=1}^n \lambda_i g_i(t - \tau(i))
  ~~~ \mbox{and} ~~~
  x(t + 1) = \argmin_{x \in \xdomain} \Big\{\<z(t+1), x\> +
  \frac{1}{\stepsize(t+1)} \prox(x) \Big\}
  \label{eqn:da-dist-update}
\end{equation}
while for mirror descent, the update is
\begin{equation}
  \lambdagrad(t)
  = \sum_{i=1}^n \lambda_i g_i(t - \tau(i))
  ~~~ \mbox{and} ~~~
  x(t + 1) = \argmin_{x \in \xdomain} \Big\{\<\lambdagrad(t), x\>
  + \frac{1}{\stepsize(t)} \divergence(x, x(t))\Big\}.
  \label{eqn:md-dist-update}
\end{equation}
The next two theorems build on the proofs of
Theorems~\ref{theorem:da-delayed-stoch} and~\ref{theorem:md-delayed-stoch},
combining several techniques. We provide the proof of
Theorem~\ref{theorem:da-distributed-delay} in
Sec.~\ref{sec:proof-da-distributed-delay}, omitting the proof of
Theorem~\ref{theorem:md-distributed-delay} as it follows in a similar way from
Theorem~\ref{theorem:md-delayed-stoch}.
\begin{theorem}
  \label{theorem:da-distributed-delay}
  Let the sequence $x(t)$ be defined by the update
  \eqref{eqn:da-dist-update}.  Under
  assumptions~\ref{assumption:lipschitz-func},~\ref{assumption:lipschitz-grad}
  and~\ref{assumption:compact}, let $\frac{1}{\stepsize(t)} = L +
  \extrastep(t)$ and $\extrastep(t) \propto \sqrt{t+\tau}$ or
  $\extrastep(t) \equiv \extrastep$ for all $t$. Then
  \begin{align*}
    \E\left[\sum_{t=1}^T f(x(t + 1)) - T f(x^*)\right]
    & \le L\radius^2 + \extrastep(T) \radius^2
    + 2 \sum_{i=1}^n \lambda_i \tau(i) G\radius
    + 2 \sum_{i=1}^n \lambda_i LG^2 (\tau(i) + 1)^2
    \sum_{t=1}^T \frac{1}{\extrastep(t-\tau)^2} \\
    & \qquad ~ +
    \sum_{t=1}^T \frac{1}{2 \extrastep(t)}
    \E \bigg\|\sum_{i=1}^n \lambda_i [\nabla f(x(t - \tau(i)))
      - g_i(t - \tau(i))]\bigg\|_*^2.
  \end{align*}
\end{theorem}
\begin{theorem}
  \label{theorem:md-distributed-delay} Use the same conditions as
  Theorem~\ref{theorem:da-distributed-delay}, but assume that $x(t)$ is
  defined by the update~\eqref{eqn:md-dist-update} and
  $\divergence(x^*, x) \le \radius^2$ for all $x \in \xdomain$. Then  
  \begin{align*}
    \E\left[\sum_{t=1}^T f(x(t + 1)) - Tf(x^*)\right]
    & \le 2 \radius^2(L + \extrastep(T)) + 2 \sum_{i=1}^n \lambda_i \tau(i)
    G\radius
    + 2 \sum_{i=1}^n \lambda_i L G^2(\tau(i) + 1)^2
    \sum_{t=1}^T \frac{1}{\extrastep(t-\tau)^2} \\
    & \qquad ~ + \sum_{t=1}^T \frac{1}{2\extrastep(t)}
    \E\bigg\|\sum_{i=1}^n \lambda_i[\nabla f(x(t - \tau(i)))
      - g_i(t - \tau(i))]
    \bigg\|_*^2.
  \end{align*}
\end{theorem}
\noindent
The consequences of Theorems~\ref{theorem:da-distributed-delay}
and~\ref{theorem:md-distributed-delay} are powerful, as we illustrate in
the next section.

\section{Distributed Optimization}
\label{sec:distributed}

We now turn to what we see as the main purpose and application of the
above results: developing robust and efficient algorithms for
distributed stochastic optimization. Our main motivations here are
machine learning and statistical applications
where the data is so large that it cannot fit on a single
computer. Examples of the form~\eqref{eqn:stochastic-objective}
include logistic regression (for background,
see~\cite{HastieTiFr01}), where the task is to learn a linear
classifier that assigns labels in $\{-1, +1\}$ to a series of
examples, in which case we have the objective
$F(x; \statsample) = \log[1 + \exp(\<\statsample, x\>)]$ as described
in Sec.~\ref{sec:setup-cent}\eqref{item:logistic-regression};
or linear regression, where $\statsample = (a, b) \in \R^d \times \R$
and $F(x; \statsample) = \half [b - \<a, x\>]^2$ as described in
Sec.~\ref{sec:setup-cent}\eqref{item:least-squares}.
Both objectives satisfy assumptions~\ref{assumption:lipschitz-func}
and~\ref{assumption:lipschitz-grad} as discussed earlier. We consider
both stochastic and online/streaming scenarios for such problems. In
the simplest setting, the distribution $P$ in the
objective~\eqref{eqn:stochastic-objective} is the empirical
distribution over an observed dataset, that is,
\begin{equation*}
  f(x) = \frac{1}{N}\sum_{i=1}^N F(x;\statsample_i).
\end{equation*}
We divide the $N$ samples among $n$ workers so that each worker has an
$N/n$-sized subset of data. 
In streaming
applications, the distribution $P$ is the unknown distribution
generating the data, and each worker receives a stream of independent
data points $\statsample \sim P$.  Worker $i$ uses its subset of the
data, or its stream, to compute $g_i$, an estimate of the gradient
$\nabla f$ of the global $f$. We make the simplifying assumption that
$g_i$ is an unbiased estimate of $\nabla f(x)$, which is satisfied,
for example, when each worker receives an
independent stream of samples or computes the gradient $g_i$ based on
samples picked at random without replacement from its subset of the
data.

The architectural assumptions we make are natural and based off of
master/worker topologies, but the convergence results in
Section~\ref{sec:main-results} allow us to give procedures robust to
delay and asynchrony. We consider two protocols: in the first,
workers compute and
communicate asynchronously and independently with the master, and in the
second, workers are at different distances from the master
and communicate with time lags proportional to their distances. We
show in the latter part of this section that the convergence rates of
each protocol, when applied in an $n$-node network, are $\order(1 /
\sqrt{nT})$ for $n$-node networks (though lower order terms are
different for each).

Before describing our architectures, we note that perhaps the simplest
master-worker scheme is to have each worker simultaneously
compute a stochastic gradient and send it to the master, which takes a
gradient step on the averaged gradient. While the $n$ gradients are
computed in parallel, accumulating and averaging $n$ gradients at the
master takes $\Omega(n)$ time, offsetting the gains of
parallelization. Thus we consider alternate architectures that are
robust to delay.

\paragraph{Cyclic Delayed Architecture}
This protocol is the delayed update algorithm mentioned in the introduction,
and it parallelizes computation of (estimates of) the gradient $\nabla
f(x)$. Formally, worker $i$ has parameter $x(t)$ and computes $g_i(t) =
F(x(t);\statsample_i(t))$, where $\statsample_i(t)$ is a random variable
sampled at worker $i$ from the distribution $P$. The master maintains a
parameter vector $x \in \xdomain$. The algorithm proceeds in rounds,
cyclically pipelining updates. The algorithm begins by initiating
gradient computations at different workers at slightly offset times. At time
$t$, the master receives gradient information at a $\tau$-step delay from some
worker, performs a parameter update, and passes the updated central parameter
$x(t + 1)$ back to the worker. Other workers do not see this update and
continue their gradient computations on stale parameter vectors. In the
simplest case, each node suffers a delay of $\tau = n$, though our earlier
analysis applies to random delays throughout the network as well. Recall
Fig.~\ref{fig:master-worker} for a graphic description of the process.

\begin{figure}[t]
  \begin{center}
    \begin{tabular}{cc}
      \includegraphics[width=.4\columnwidth]{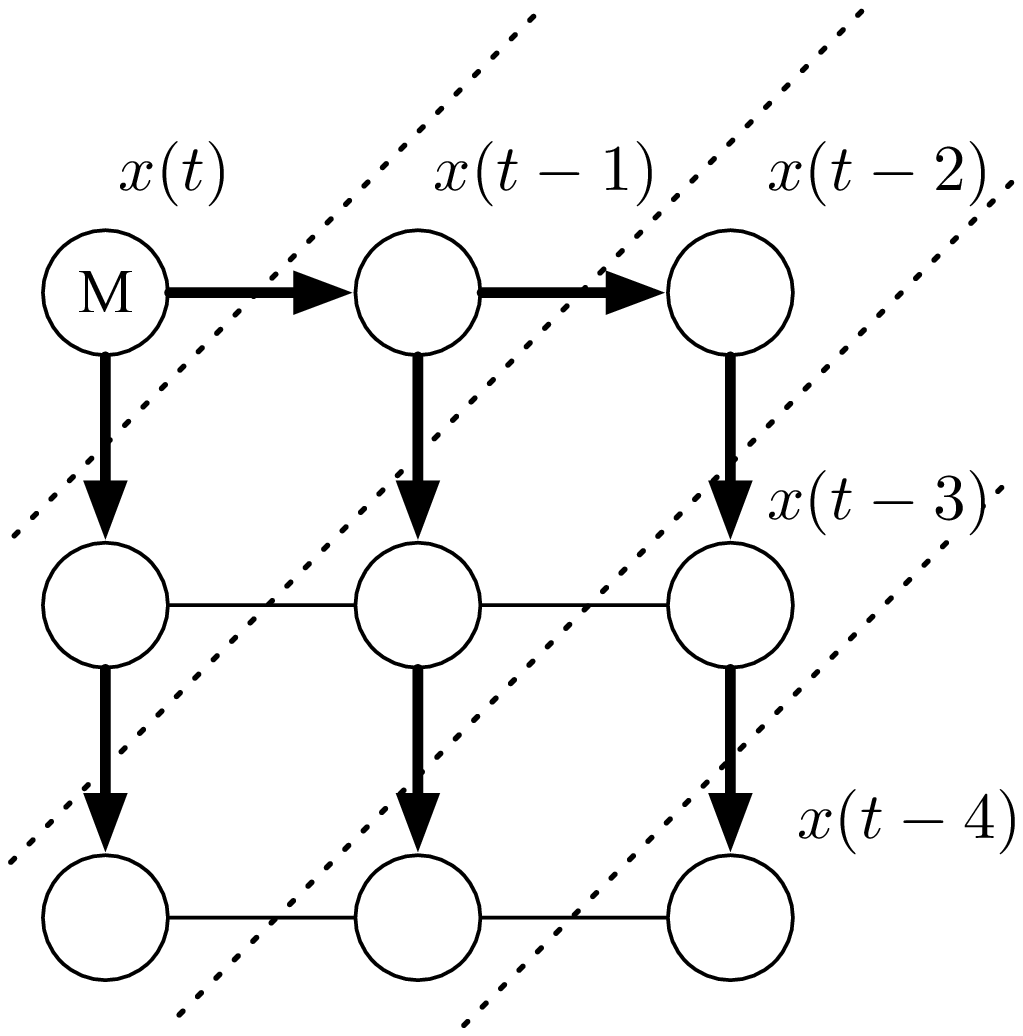} &
      \includegraphics[width=.4\columnwidth]{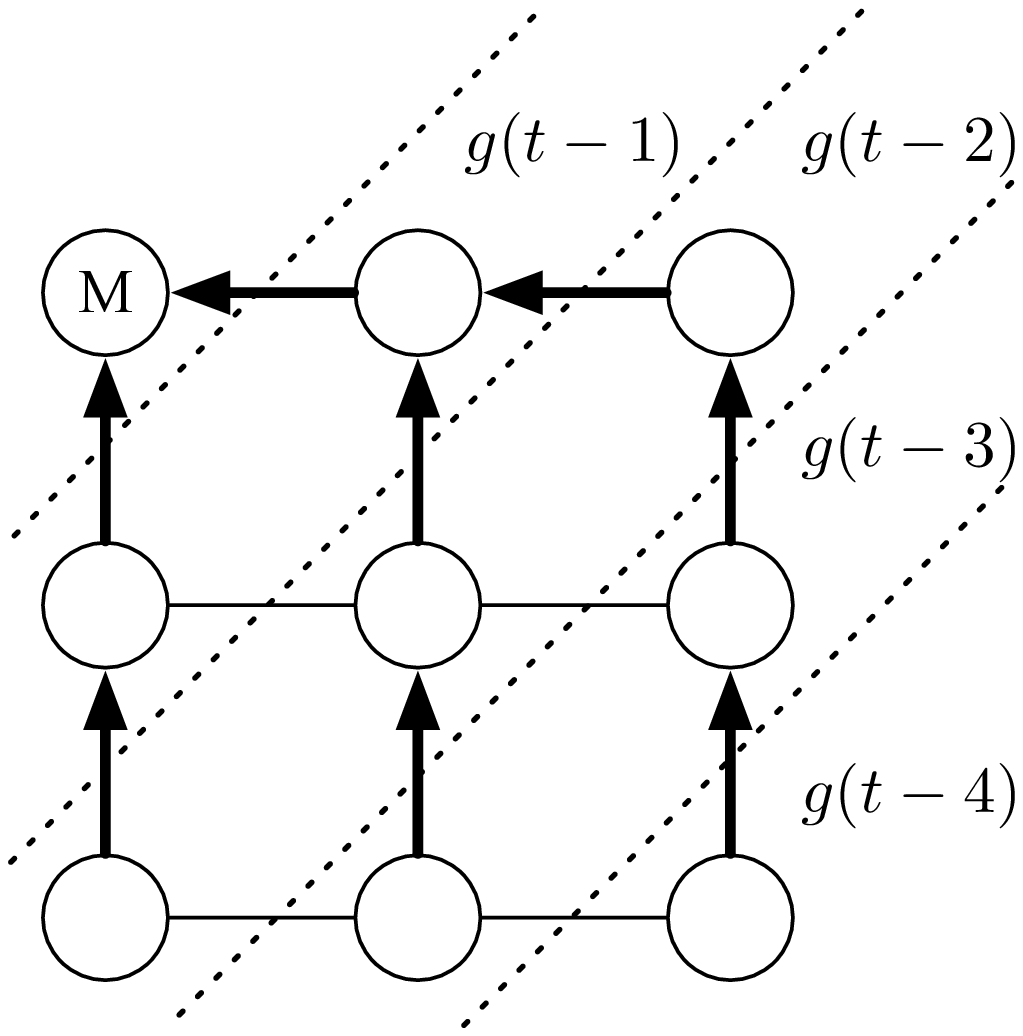} \\
      (a) & (b)
    \end{tabular}
  \end{center}
  \caption{\label{fig:averaging-network} Master-worker averaging network.
    (a): parameters stored at different distances from master node at
    time $t$. A node at distance $d$ from master has the parameter
    $x(t-d)$.  (b): gradients computed at different nodes. A node at
    distance $d$ from master computes gradient $g(t-d)$.}
\end{figure}

\begin{figure}[t]
  \begin{center}
    \includegraphics[width=.6\columnwidth]{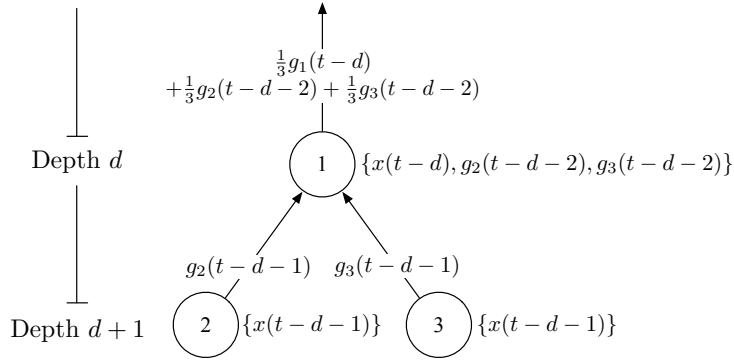}
    \caption{\label{fig:node-averaging} Communication of gradient information
      toward master node at time $t$ from node 1 at distance $d$ from
      master. Information stored at time $t$ by node $i$ in brackets to right
      of node.}
  \end{center}
\end{figure}

\paragraph{Locally Averaged Delayed Architecture}
At a high level, the protocol we now describe combines the delayed updates of
the cyclic delayed architecture with averaging techniques of previous
work~\cite{NedicOz09,DuchiAgWa10}. We assume a network $\graph = (\vertices,
\edges)$, where $\vertices$ is a set of $n$ nodes (workers) and $\edges$ are
the edges between the nodes. We select one of the nodes as the master, which
maintains the parameter vector $x(t) \in \xdomain$ over time.

The algorithm works via a series of multicasting and aggregation steps on a
spanning tree rooted at the master node. In the first phase, the algorithm
broadcasts from the root towards the leaves. At step $t$ the master sends its
current parameter vector $x(t)$ to its immediate neighbors. Simultaneously,
every other node broadcasts its current parameter vector (which, for a depth
$d$ node, is $x(t-d)$) to its children in the spanning tree. See
Fig.~\ref{fig:averaging-network}(a). Every worker receives its new parameter
and computes its local gradient at this parameter. The second part of the
communication in a given iteration proceeds from leaves toward the root. The
leaf nodes communicate their gradients to their parents. The parent takes the
gradients of the leaf nodes from the previous round (received at iteration $t
- 1$) and averages them with its own gradient, passing this averaged gradient
back up the tree. Again simultaneously, each node takes the averaged gradient
vectors of its children from the previous rounds, averages them with its
current gradient vector, and passes the result up the
spanning tree. See Fig.~\ref{fig:averaging-network}(b) and
Fig.~\ref{fig:node-averaging} for a visual description.

Slightly more formally, associated with each node $i \in \vertices$ is a delay
$\tau(i)$, which is (generally) twice its distance from the master. Fix an
iteration $t$. Each node $i \in \vertices$ has an out of date parameter vector
$x(t - \tau(i)/2)$, which it sends further down the tree to its children. So,
for example, the master node sends the vector $x(t)$ to its children, which
send the parameter vector $x(t - 1)$ to their children, which in turn send
$x(t - 2)$ to their children, and so on. Each node computes
\begin{equation*}
  g_i(t - \tau(i)/2) = \nabla F(x(t - \tau(i)/2); \statsample_i(t)),
\end{equation*}
where $\statsample_i(t)$ is a random variable sampled at node $i$ from the
distribution $P$. The communication back up the hierarchy proceeds as follows:
the leaf nodes in the tree (say at depth $d$) send the gradient vectors $g_i(t
- d)$ to their immediate parents in the tree. At the previous iteration $t -
1$, the parent nodes received $g_i(t - d - 1)$ from their children, which they
average with their own gradients $g_i(t - d + 1)$ and pass to their parents,
and so on. The master node at the root of the tree receives an average
of delayed gradients from the entire tree, with each gradient having a
potentially different delay, giving rise to updates of the
form~(\ref{eqn:da-dist-update}) or~(\ref{eqn:md-dist-update}).

\subsection{Convergence rates for delayed distributed minimization}
\label{sec:delayeddist}

Having described our architectures, we can now give corollaries to the
theoretical results from the previous sections that show it is possible
to achieve asymptotically faster rates (over centralized procedures) using
distributed algorithms even without imposing synchronization requirements.  We
allow workers to pipeline updates by computing asynchronously and in parallel,
so each worker can compute low variance estimate of the gradient $\nabla
f(x)$.

We begin with a simple corollary to the results in
Sec.~\ref{sec:delayedcyclic}. We
ignore the constants $L$, $G$, $\radius$, and $\stddev$, which are not
dependent on the characteristics of the network. We also assume that
each worker uses $m$ independent samples of $\statsample \sim P$ to
compute the stochastic gradient as
\begin{align*}
  g_i(t) = \frac{1}{m}\sum_{j=1}^m\nabla F(x(t);\statsample_i(j)).
\end{align*}
 Using the cyclic protocol as in
Fig.~\ref{fig:master-worker},
Theorems~\ref{theorem:da-delayed-stoch} and~\ref{theorem:md-delayed-stoch}
give the following result.
\begin{corollary}
  \label{corollary:cyclic-architecture}
  Let $\prox(x) = \half \ltwo{x}^2$, assume the conditions in
  Corollary~\ref{corollary:random-tau}, and assume that each worker
  uses $m$ samples $\statsample \sim P$ to compute the gradient it
  communicates to the master. Then with the choice $\extrastep(t) =
  \sqrt{T} / \sqrt{m}$ either of the
  updates~\eqref{eqn:da-delayed-update}
  or~\eqref{eqn:md-delayed-update} satisfy
  \begin{equation*}
    \E[f(\what{x}(T))] - f(x^*) = \order\bigg(\frac{B^2}{T} +
    \frac{1}{\sqrt{Tm}} + \frac{B^2 m}{T}\bigg).
  \end{equation*}
\end{corollary}
\begin{proof}
  The corollary follows straightforwardly from the realization that the
  variance $\stddev^2 = \E[\ltwo{\nabla f(x) - g_i(t)}^2] =
  \E[\ltwo{\nabla f(x) - \nabla F(x;\statsample)}^2]/m = \order(1/m)$
  when workers use $m$ independent stochastic gradient samples. 
\end{proof}
\noindent
In the above corollary, so long as the bound on the delay $B$
satisfies, say, $B = o(T^{1/4})$, then the last term in the bound is
asymptotically negligible, and we achieve a convergence rate of
$\order(1 / \sqrt{Tm})$.

The cyclic delayed architecture has the drawback that information from a
worker can take $\order(n)$ time to reach the master. While the algorithm is
robust to delay and does not need lock-step coordination of workers, the
downside of the architecture is that the essentially $n^2 m/ T$ term in the
bounds above can be quite large. Indeed, if each worker computes its gradient
over $m$ samples with $m \approx n$---say to avoid idling of workers---then
the cyclic architecture has convergence rate $\order(n^3 / T +
1/\sqrt{nT})$.
For moderate $T$ or large
$n$, the delay penalty $n^3/T$ may dominate $1 / \sqrt{nT}$,
offsetting the gains of parallelization.


To address the large $n$ drawback, we turn our attention to the
locally averaged architecture described by
Figs.~\ref{fig:averaging-network} and~\ref{fig:node-averaging}, where
delays can be smaller since they depend only on the height of a
spanning tree in the network. The algorithm requires more
synchronization than the cyclic architecture but still performs
limited local communication. Each worker computes $g_i(t -
\tau(i)) = \nabla F(x(t-\tau(i)); \statsample_i(t))$ where $\tau(i)$
is the delay of worker $i$ from the master and $\statsample_i \sim
P$. As a result of the communication procedure, the master receives a
convex combination of the stochastic gradients evaluated at each
worker $i$, for which we gave results in
Section~\ref{sec:combo-delays}.

In this architecture, the master receives gradients of the
form $g_\lambda(t) = \sum_{i=1}^n \lambda_i g_i(t - \tau(i))$ for some
$\lambda$ in the simplex, which puts us in the setting of
Theorems~\ref{theorem:da-distributed-delay}
and~\ref{theorem:md-distributed-delay}. We now make the reasonable
assumption that the gradient errors $\nabla f(x(t)) - g_i(t)$ are
uncorrelated across the nodes in the network.\footnote{Similar results
  continue to hold under weak correlation.} In statistical
applications, for example, each worker may own independent data or
receive streaming data from independent sources; more generally, each
worker can simply receive independent samples $\statsample_i\sim
P$. We also set $\prox(x) = \half \ltwo{x}^2$, and observe
\begin{equation*}
  \E \bigg\|\sum_{i=1}^n \lambda_i \nabla f(x(t - \tau(i))) - g_i(t -
  \tau(i))\bigg\|_2^2 = \sum_{i=1}^n \lambda_i^2 \E \ltwo{\nabla f(x(t
    - \tau(i))) - g_i(t - \tau(i))}^2.
\end{equation*}
This gives the following corollary to
Theorems~\ref{theorem:da-distributed-delay}
and~\ref{theorem:md-distributed-delay}.
\begin{corollary}
  \label{corollary:distributed-delay}
  Set $\lambda_i = \ninv$ for all $i$, $\prox(x) = \half \ltwo{x}^2$,
  and $\extrastep(t) = \sigma \sqrt{t+\tau}/ \radius \sqrt{n}$.  Let
  $\bar{\tau}$ and $\overline{\tau^2}$ denote the average of the
  delays $\tau(i)$ and $\tau(i)^2$, respectively. Under the conditions
  of Theorem~\ref{theorem:da-distributed-delay}
  or~\ref{theorem:md-distributed-delay},
  \begin{equation*}
    \E \left[\sum_{t=1}^T f(x(t + 1)) - T f(x^*)\right] =
    \order\left(L\radius^2 + \bar{\tau} G\radius + \frac{L G^2
      \radius^2 n\overline{\tau^2}}{\sigma^2} \log T +
    \frac{\radius\sigma}{\sqrt{n}} \sqrt{T}\right).
  \end{equation*}
\end{corollary}
\noindent
The $\log T$ multiplier can be reduced to a constant if we set
$\extrastep(t) \equiv \stddev \sqrt{T} / R \sqrt{n}$. By using the
averaged sequence $\what{x}(T)$~\eqref{eqn:avg-x}, Jensen's inequality
gives that asymptotically $\E[f(\what{x}(T))] - f(x^*) = \order(1 /
\sqrt{Tn})$, which is an optimal dependence on the number of samples
$\statsample$ calculated by the method. We also observe in this
architecture, the delay $\tau$ is bounded by the graph diameter
$\diam$, giving us the bound:
\begin{equation}
  \E \left[\sum_{t=1}^T f(x(t + 1)) - T f(x^*)\right] =
  \order\left(L\radius^2 + \diam G\radius + \frac{L G^2 \radius^2
    n\diam^2}{\sigma^2} \log T + \frac{\radius\sigma}{\sqrt{n}}
  \sqrt{T}\right).
  \label{eqn:delayed-diam}
\end{equation}

The above corollaries are general and hold irrespective of the relative costs
of communication and computation. However, with knowledge of the costs, we can
adapt the stepsizes slightly to give better rates of convergence when $n$ is
large or communication to the master node is expensive.  For now, we focus on
the cyclic architecture (the setting of
Corollary~\ref{corollary:cyclic-architecture}), though the same principles
apply to the local averaging scheme. Let $\comcost$ denote the cost
of communicating between the master and workers in terms of the time to
compute a single gradient sample, and assume that we set $m = \comcost n$, so
that no worker node has idle time. For simplicity, we let the delay be
non-random, so $B = \tau = n$.  Consider the choice $\extrastep(t) =
\extrastep\sqrt{T/(\comcost n)}$ for the damping stepsizes, where $\extrastep
\geq 1$. This setting in Theorem~\ref{theorem:da-delayed-stoch} gives
\begin{equation*}
  \E[f(\what{x}(T))] - f(x^*) = \order\left(\frac{\extrastep^2
    \comcost n^3}{T} + \frac{\extrastep}{\sqrt{T \comcost n}} +
  \frac{1}{\extrastep\sqrt{T \comcost n}}\right) =
  \order\left(\frac{\extrastep^2 \comcost n^3}{T} +
  \frac{\extrastep}{\sqrt{T \comcost n}}\right),
\end{equation*}
where the last equality follows since $\extrastep \geq 1$. Optimizing for
$\extrastep$ on the right yields
\begin{equation}
  \extrastep = \max\left\{\frac{n^{7/6}}{T^{1/6} \comcost^{1/2}},
  1\right\} ~~~\mbox{and}~~~ \E[f(\what{x}(T))] - f(x^*) =
  \order\left(\min\left\{\frac{n^{2/3}}{T^{2/3}},\frac{n^3}{T}\right\}
  + \frac{1}{\sqrt{T \comcost n}}\right).
  \label{eqn:smaller-cyclic-bound}
\end{equation}
The convergence rates thus follow two regimes. When $T \leq n^7 /
\comcost^3$, we have convergence rate $\order(n^{2/3}/T^{2/3})$, while
once $T > n^7 / \comcost^3$, we attain $\order(1/\sqrt{T\comcost n})$
convergence. Roughly, in time proportional to $T \comcost$, we achieve
optimization error $1 / \sqrt{T \comcost n}$, which is order-optimal
given that we can compute a total of $T \comcost n$ stochastic
gradients~\cite{AgarwalBaRaWa10}. The scaling of this bound is nicer
than that previously: the dependence on network size is at worst
$n^{2/3}$, which we obtain by increasing the damping factor
$\extrastep(t)$---and hence decreasing the stepsize $\stepsize(t) = 1
/ (L + \extrastep(t))$---relative to the setting of
Corollary~\ref{corollary:cyclic-architecture}.
We remark that applying the same technique to
Corollary~\ref{corollary:distributed-delay} gives convergence rate
scaling as the smaller of 
$\order((\diam/T)^{2/3} + 1/\sqrt{T \comcost n})$ and
$\order((nC\diam/T + 1/\sqrt{T \comcost n})$. Since the diameter
$\diam \le n$, this is faster than the cyclic architecture's
bound~\eqref{eqn:smaller-cyclic-bound}.

\subsection{Running-time comparisons}
\label{sec:wall-clock}

Having derived the rates of convergence of the different distributed
procedures above, we now explicitly study the running times of the
centralized stochastic gradient algorithms~\eqref{eqn:da-update}
and~\eqref{eqn:md-update}, the cyclic delayed protocol with the
updates~\eqref{eqn:da-delayed-update}
and~\eqref{eqn:md-delayed-update}, and the locally averaged
architecture with the updates~\eqref{eqn:da-dist-update}
and~\eqref{eqn:md-dist-update}.
To make comparisons more cleanly, we avoid constants, assuming without
loss that the variance bound $\sigma^2$ on $\E \norm{\nabla f(x) -
  \nabla F(x; \statsample)}^2$ is $1$, and that sampling $\statsample
\sim P$ and evaluating $\nabla F(x; \statsample)$ requires one unit of
time. Noting that $\E[\nabla F(x; \statsample)] = \nabla f(x)$, it is
clear that if we receive $m$ uncorrelated samples of $\statsample$,
the variance $\E\sltwo{\nabla f(x) - \frac{1}{m} \sum_{j=1}^m \nabla
  F(x; \statsample_j)}^2 \le \frac{1}{m}$.

Now we state our assumptions on the relative times used by each
algorithm. Let $T$ be the number of units of time allocated to each
algorithm, and let the centralized, cyclic delayed and locally
averaged delayed algorithms complete $\tcent$, $\tcycle$ and $\tdist$
iterations, respectively, in time $T$. It is clear that $\tcent =
T$. We assume that the distributed methods use $\mcycle$ and $\mdist$
samples of $\statsample \sim P$ to compute stochastic gradients and
that the delay $\tau$ of the cyclic algorithm is $n$.
For concreteness, we assume that communication is of the same order as
computing the gradient of one sample $\nabla F(x; \statsample)$ so
that $\comcost = 1$. In the cyclic setup of
Sec.~\ref{sec:delayedcyclic}, it is reasonable to assume that $\mcycle
= \Omega(n)$ to avoid idling of workers
(Theorems~\ref{theorem:da-delayed-stoch}
and~\ref{theorem:md-delayed-stoch}, as well as the
bound~\eqref{eqn:smaller-cyclic-bound}, show it is asymptotically
beneficial to have $\mcycle$ larger, since $\stddevcycle^2 = 1 /
\mcycle$). For $\mcycle = \Omega(n)$, the master requires
$\frac{\mcycle}{n}$ units of time to receive one gradient update, so
$\frac{\mcycle}{n} \tcycle = T$. In the locally delayed framework, if
each node uses $\mdist$ samples to compute a gradient, the master
receives a gradient every $\mdist$ units of time, and hence
$\mdist\tdist = T$. Further, $\stddevdist^2 = 1 / \mdist$.
We summarize our assumptions by saying that in $T$
units of time, each algorithm performs the following number of
iterations:
\begin{equation}
  \label{eqn:updates-per-time}
  \tcent = T, ~~~~ \tcycle = \frac{Tn}{\mcycle}, ~~~~ {\rm and} ~~~~
  \tdist = \frac{T}{\mdist}.
\end{equation}


\begin{table}[t]
  \begin{center}
    \begin{tabular}{|c|c|}
      \hline Centralized (\ref{eqn:da-update}, \ref{eqn:md-update}) &
      $\displaystyle{ \E f(\xavg) - f(x^*) = \order\bigg(
        \sqrt{\frac{1}{T}}\bigg)}$ \\ \hline Cyclic
      (\ref{eqn:da-delayed-update}, \ref{eqn:md-delayed-update}) &
      $\displaystyle{
        \order\left(\min\left(\frac{n^{2/3}}{T^{2/3}},\frac{n^3}{T}\right)
        + \frac{1}{\sqrt{T n}}\right)}$\\
      \hline Local (\ref{eqn:da-dist-update},
      \ref{eqn:md-dist-update}) & $\displaystyle{ \E f(\xavg) - f(x^*)
        = \order\left(
        \min\left(\frac{\diam^{2/3}}{T^{2/3}},\frac{n
          \overline{\tau^2}}{T}\right) + \frac{1}{\sqrt{n
            T}}\right)}$ \\ \hline
    \end{tabular}
    \caption{\label{table:wall-clock} Upper bounds on $\E f(\xavg) - f(x^*)$
      for three computational architectures, where $\xavg$ is the output of
      each algorithm after $T$ units of time. Each algorithm runs for the
      amount of time it takes a centralized stochastic algorithm to perform
      $T$ iterations as in~\eqref{eqn:updates-per-time}. Here $\diam$ is the
      diameter of the network, $n$ is the number of nodes, and
      $\overline{\tau^2} = \frac{1}{n} \sum_{i=1}^n \tau(i)^2$ is the average
      squared communication delay for the local averaging architecture.
      Bounds for the cyclic architecture assume delay $\tau
      = n$.}
  \end{center}
\end{table}

Plugging the above iteration counts into the earlier bound
\eqref{eqn:smooth-da-stochastic} and
Corollaries~\ref{corollary:cyclic-architecture}
and~\ref{corollary:distributed-delay} via the sharper
result~\eqref{eqn:smaller-cyclic-bound}, we can provide upper bounds (to
constant factors) on the expected optimization accuracy after $T$ units of
time for each of the distributed architectures as in
Table~\ref{table:wall-clock}.
Asymptotically in the number of units of time $T$, both the cyclic and locally
communicating stochastic optimization schemes have the same convergence
rate. However, topological considerations show that the locally communicating
method (Figs.~\ref{fig:averaging-network} and~\ref{fig:node-averaging}) has
better performance than the cyclic architecture, though it requires more
worker coordination. Since
the lower order terms matter only for large $n$ or small $T$, we compare
the terms $n^{2/3}/T^{2/3}$ and $\diam^{2/3}/T^{2/3}$ for the cyclic
and locally averaged algorithms, respectively. Since $\diam \le n$ for
any network, the locally averaged algorithm always guarantees better
performance than the cyclic algorithm. For
specific graph topologies, however, we can quantify the
time improvements:
\begin{itemize}
\item $n$-node cycle or path: $\diam = n$ so that both methods have
  the same convergence rate.
\item $\sqrt{n}$-by-$\sqrt{n}$ grid: $\diam = \sqrt{n}$, so the
  distributed method has a factor of $n^{2/3}/n^{1/3} = n^{1/3}$
  improvement over the cyclic architecture.
\item Balanced trees and expander graphs: $\diam = \order(\log n)$, so
  the distributed method has a factor---ignoring logarithmic
  terms---of $n^{2/3}$ improvement over cyclic.
\end{itemize}

Naturally, it is possible to modify our assumptions. In a network in
which communication is cheap, or conversely, in a problem for which
the computation of $\nabla F(x; \statsample)$ is more expensive than
communication, then the number of samples $\statsample \sim P$ for
which which each worker computes gradients is small. Such problems are
frequent in statistical machine learning, such as when learning
conditional random field models, which are useful in natural language
processing, computational biology, and other application
areas~\cite{LaffertyMcPe01}. In this case, it is reasonable to have
$\mcycle = \order(1)$, in which case $\tcycle = Tn$ and
the cyclic delayed architecture
has stronger convergence guarantees of $\order(\min\{n^2/T, 1/T^{2/3}\} + 1 /
\sqrt{Tn})$.  In any case, both non-centralized protocols enjoy
significant asymptotically faster convergence rates for stochastic
optimization problems in spite of asynchronous delays.

\section{Numerical Results}
\label{sec:experiments}

Though this paper focuses mostly on the theoretical analysis of the methods we
have presented, it is important to understand the practical aspects of the
above methods in solving real-world tasks and problems with real data. To that
end, we use the cyclic delayed method~\eqref{eqn:da-dist-update} to
solve a common statistical machine learning problem. Specifically, we focus
on solving the logistic regression problem
\begin{equation}
  \label{eqn:logreg-example}
  \min_x ~ f(x) = \frac{1}{N} \sum_{i=1}^N \log(1 + \exp(-b_i\<a_i, x\>))
  ~~~ \mbox{subject~to~} \ltwo{x} \le \radius.
\end{equation}
We use the Reuters RCV1 dataset~\cite{LewisYaRoLi04}, which consists of $N
\approx 800000$ news articles, each labeled with some combination of the four
labels economics, government, commerce, and medicine. In the above example,
the vectors $a_i \in \{0, 1\}^d$, $d \approx 10^5$, are feature vectors
representing the words in each article, and the labels $b_i$ are $1$ if the
article is about government, $-1$ otherwise.
\begin{figure}[t]
  \begin{center}
    \includegraphics[width=.45\columnwidth]{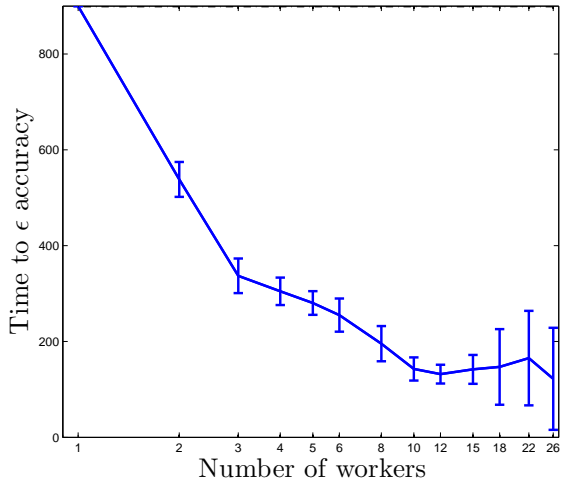}
    \caption{\label{fig:single-example-communication} Optimization performance
      of the delayed cyclic method~\eqref{eqn:da-delayed-update} for the
      Reuters RCV1 dataset when we assume that the cost of communication to
      the master is the same as computing the gradient of one term in the
      objective~\eqref{eqn:logreg-example}. The number of samples $m$ computed
      is equal to $n$ for each worker. Plotted is the estimated time to
      $\epsilon$-accuracy as a function of number of workers $n$.}
  \end{center}
\end{figure}

We simulate the cyclic delayed optimization
algorithm~\eqref{eqn:da-delayed-update} for the
problem~\eqref{eqn:logreg-example} for several choices of the number of
workers $n$ and the number of samples $m$ computed at each worker.  We
summarize the results of our experiments in
Figure~\ref{fig:single-example-communication}.
To generate the figure, we fix an
$\epsilon$ (in this case, $\epsilon = .05$), then measure the time it takes
the stochastic algorithm~\eqref{eqn:da-delayed-update} to output an $\what{x}$
such that $f(\what{x}) \le \inf_{x \in \xdomain} f(x) + \epsilon$. We perform
each experiment ten times.

After computing the number of iterations required to achieve
$\epsilon$-accuracy, we convert the results to running time by assuming it
takes one unit of time to compute the gradient of one term in the sum defining
the objective~\eqref{eqn:logreg-example}. We also assume that it takes $1$
unit of time, i.e.\ $\comcost = 1$, to communicate from one of the workers to
the master, for the master to perform an update, and communicate back to one
of the workers. In an $n$ node system where each worker computes $m$ samples
of the gradient, the master receives an update every $\max\{\frac{m}{n}, 1\}$
time units. A centralized algorithm computing $m$ samples of its gradient
performs an update every $m$ time units. By multiplying the number of
iterations to $\epsilon$-optimality by $\max\{\frac{m}{n}, 1\}$ for the
distributed method and by $m$ for the centralized, we can estimate the amount
of time it takes each algorithm to achieve an $\epsilon$-accurate solution.

We now turn to discussing Figure~\ref{fig:single-example-communication}.  The
delayed update~\eqref{eqn:da-delayed-update} enjoys speedup (the ratio of time
to $\epsilon$-accuracy for an $n$-node system versus the centralized
procedure) nearly linear in the number $n$ of worker machines until $n \ge 15$
or so. Since we use the stepsize choice $\extrastep(t) \propto \sqrt{t
  / n}$, which yields the predicted convergence rate given by
Corollary~\ref{corollary:cyclic-architecture}, the $n^2 m / T \approx n^3 / T$
term in the convergence rate presumably becomes non-negligible for larger
$n$. This expands on earlier experimental work with a similar
method~\cite{LangfordSmZi09}, which experimentally demonstrated linear speedup
for small values of $n$, but did not investigate larger network sizes.
Roughly, as predicted by our theory, for non-asymptotic regimes the
cost of communication and delays due to using $n$ nodes mitigate some of the
benefits of parallelization. Nevertheless, as our analysis shows,
allowing delayed and asynchronous updates still gives significant performance
improvements.

\section{Delayed Updates for Smooth Optimization}
\label{sec:smooth-delay}

In this section, we prove Theorems~\ref{theorem:da-delayed-stoch}
and~\ref{theorem:md-delayed-stoch}. We collect in
Appendix~\ref{sec:technical-prox} a few technical results relevant to our
proof; we will refer to results therein without comment.
Before proving either theorem, we state the lemma that is the key to our
argument. Lemma~\ref{lemma:bound-error-smooth} shows that certain
gradient-differencing terms are essentially of second order. As a consequence,
when we combine the results of the lemma with
Lemma~\ref{lemma:xt-difference-bound}, which bounds $\E[\norm{x(t) - x(t +
    \tau)}^2]$, the gradient differencing terms become $\order(\log T)$ for
step size choice $\extrastep(t) \propto \sqrt{t}$, or $\order(1)$ for
$\extrastep(t) \equiv \extrastep \sqrt{T}$.
\begin{lemma}
  \label{lemma:bound-error-smooth}
  Let assumptions~\ref{assumption:lipschitz-func}
  and~\ref{assumption:lipschitz-grad} on the function $f$ and the
  compactness assumption~\ref{assumption:compact} hold. Then for any
  sequence $x(t)$
  \begin{equation*}
    \sum_{t=1}^T \<\nabla f(x(t)) - \nabla f(x(t - \tau)), x(t + 1) - x^*\>
    \le \frac{L}{2}\sum_{t=1}^T\norm{x(t - \tau) - x(t+1)}^2
    + 2 \tau G\radius.
  \end{equation*}
  Consequently, if $\E[\norm{x(t) - x(t+1)}^2] \leq \kappa(t)^2G^2$ for a
  non-increasing sequence $\kappa(t)$,
  \begin{equation*}
    \E\bigg[\sum_{t=1}^T
      \<\nabla f(x(t)) - \nabla f(x(t - \tau)), x(t + 1) - x^*\>\bigg]
    \le \frac{LG^2 (\tau + 1)^2}{2} \sum_{t=1}^T \kappa(t-\tau)^2
    + 2 \tau G\radius.
  \end{equation*}
\end{lemma}
\begin{proof} 
  The proof follows by using a few Bregman divergence identities to rewrite
  the left hand side of the above equations, then recognizing that the result
  is close to a telescoping sum. Recalling the definition of a Bregman
  divergence~\eqref{eqn:def-bregman-divergence}, we note the following
  well-known four term equality, a consequence of straightforward algebra:
  for any $a, b, c, d$,
  \begin{equation}
    \label{eqn:three-term-bregman}
    \<\nabla f(a) - \nabla f(b), c - d\>
    = \divergencef(d, a) - \divergencef(d, b) - \divergencef(c, a)+
    \divergencef(c,b).  
  \end{equation}
  Using the equality \eqref{eqn:three-term-bregman}, we see that 
  \begin{align}
    \lefteqn{\<\nabla f(x(t)) - \nabla f(x(t - \tau)), x(t + 1) - x^*\>}
    \nonumber \\
    & = \divergencef(x^*, x(t)) - \divergencef(x^*, x(t - \tau))
    - \divergencef(x(t + 1), x(t)) + \divergencef(x(t + 1), x(t - \tau)).
    \label{eqn:bregmans-to-telescope}
  \end{align}
  To make \eqref{eqn:bregmans-to-telescope} useful, we note that
  the Lipschitz continuity of $\nabla f$ implies
  \begin{equation*}
    f(x(t + 1)) \le f(x(t - \tau)) + \<\nabla f(x(t - \tau)), x(t + 1)
    - x(t - \tau)\> + \frac{L}{2} \norm{x(t - \tau) - x(t + 1)}^2
  \end{equation*}
  so that recalling the definition of
  $\divergencef$~(\ref{eqn:def-bregman-divergence}) we have
  \begin{align*}
    \divergencef(x(t + 1), x(t - \tau)) & \le \frac{L}{2} \norm{x(t
      - \tau) - x(t + 1)}^2. 
  \end{align*}
  In particular, using the non-negativity of $\divergencef(x, y)$, we can
  replace \eqref{eqn:bregmans-to-telescope} with the bound
  \begin{equation*}
    \<\nabla f(x(t)) - \nabla f(x(t - \tau)), x(t + 1) - x^*\>
    \le \divergencef(x^*, x(t)) - \divergencef(x^*, x(t - \tau))
    + \frac{L}{2} \norm{x(t - \tau) - x(t + 1)}^2.
  \end{equation*}
  Summing the inequality, we see that
  \begin{equation}
    \sum_{t=1}^T \<\nabla f(x(t)) - \nabla f(x(t - \tau)),
    x(t + 1) - x^*\>
    \le \sum_{t = T - \tau + 1}^T \divergencef(x^*, x(t))
    + \frac{L}{2} \sum_{t=1}^T \norm{x(t - \tau) - x(t + 1)}^2.
    \label{eqn:second-order-sum}
  \end{equation}
  
  To bound the first Bregman divergence term, we recall that by
  Assumption~\ref{assumption:compact} and the strong convexity of $\prox$,
  $\norm{x^* - x(t)}^2 \leq 2\divergence(x^*, x(t)) \leq 2R^2$, and 
  hence the optimality of $x^*$ implies
  \begin{equation*}
    \divergencef(x^*, x(t)) = f(x^*) - f(x(t))
    - \<\nabla f(x(t)), x^* - x(t)\>
    \le \dnorm{\nabla f(x(t))} \norm{x^* - x(t)}
    \le 2 G\radius.
  \end{equation*}
  This gives the first bound of the lemma. For the second bound, using
  convexity, we see that 
  \begin{equation*}
    \norm{x(t - \tau) - x(t + 1)}^2
    \le (\tau + 1)^2 \sum_{s = 0}^\tau \frac{1}{\tau + 1}\norm{x(t - s)
      - x(t - s + 1)}^2,
  \end{equation*}
  so by taking expectations we have $\E[\norm{x(t) - x(t + \tau + 1)}^2] \le
  (\tau + 1)^2 \kappa(t-\tau)^2 G^2$. Since $\kappa$ is non-increasing 
  (by the
  definition of the update scheme) we see that the
  sum~\eqref{eqn:second-order-sum} is further bounded by $2 \tau G\radius +
  \frac{L}{2} \sum_{t=1}^T G^2 (\tau + 1)^2 \kappa(t-\tau)^2$ as desired.
\end{proof}

\subsection{Proof of Theorem~\ref{theorem:da-delayed-stoch}}
\label{sec:stoch-delayed}

The essential idea in this proof is to use convexity and smoothness
to bound $f(x(t)) - f(x^*)$, then use the sequence $\{\extrastep(t)\}$,
which decreases the stepsize $\stepsize(t)$,
to cancel variance terms. To begin, we define the error $\error(t)$
\begin{equation*}
  \error(t) \defeq \nabla f(x(t)) - g(t - \tau)
\end{equation*}
where $g(t - \tau) = \nabla F(x(t - \tau); \statsample(t)$ for some
$\statsample(t) \sim P$. Note that $\error(t)$ does not have zero
expectation, as there is a time delay. 

By using the convexity of $f$ and then the $L$-Lipschitz continuity of
$\nabla f$, for any $x^* \in \xdomain$, we have
\begin{align*}
  f(x(t)) - f(x^*)
  & \le \<\nabla f(x(t)), x(t) - x^*\> = \<\nabla f(x(t)), x(t + 1) - x^*\>
  + \<\nabla f(x(t)), x(t) - x(t + 1)\> \\
  & \le \<\nabla f(x(t)), x(t + 1) - x^*\> + f(x(t)) - f(x(t + 1))
  + \frac{L}{2} \norm{x(t) - x(t + 1)}^2,
\end{align*}
so that
\begin{align*}
  &f(x(t + 1)) - f(x^*) \le \<\nabla f(x(t)), x(t + 1) - x^*\>
  + \frac{L}{2} \norm{x(t) - x(t + 1)}^2 \\
  & = \<g(t - \tau), x(t + 1) - x^*\>
  + \<\error(t), x(t + 1) - x^*\> + \frac{L}{2} \norm{x(t) - x(t + 1)}^2 \\
  & = \<z(t+1), x(t + 1) - x^*\> - \<z(t), x(t + 1) - x^*\>
  + \<\error(t), x(t + 1) - x^*\> + \frac{L}{2} \norm{x(t) - x(t + 1)}^2.
\end{align*}
Now, by applying Lemma~\ref{lemma:solution-convexity} in
Appendix~\ref{sec:technical-prox} and the definition of the
update~\eqref{eqn:da-delayed-update}, we see that
\begin{align*}
  -\<z(t), x(t + 1) - x^*\>
  & \le - \<z(t), x(t) - x^*\>
  + \frac{1}{\stepsize(t)}\left[\prox(x(t + 1)) - \prox(x(t))\right]
  - \frac{1}{\stepsize(t)} \divergence(x(t + 1), x(t)),
\end{align*}
which implies
\begin{align}
  \lefteqn{f(x(t + 1)) - f(x^*)} \nonumber \\
  & \le \<z(t+1), x(t + 1) - x^*\>
  - \<z(t), x(t) - x^*\> + \frac{1}{\stepsize(t)}
  [\prox (x(t + 1)) - \prox(x(t))] \nonumber \\
  & ~~~~~
  - L \divergence(x(t + 1), x(t))
  - \extrastep(t) \divergence(x(t + 1), x(t))
  + \frac{L}{2} \norm{x(t) - x(t + 1)}^2
  + \<\error(t), x(t + 1) - x^*\> \nonumber \\
  & \le \<z(t+1), x(t + 1) - x^*\>
  - \<z(t), x(t) - x^*\> + \frac{1}{\stepsize(t)}
  [\prox (x(t + 1)) - \prox(x(t))] \nonumber \\
  & ~~~~ - \extrastep(t) \divergence(x(t + 1), x(t))
  + \<\error(t), x(t + 1) - x^*\>.
  \label{eqn:smooth-single-step-delayed-bound}
\end{align}
To get the bound~\eqref{eqn:smooth-single-step-delayed-bound}, we
substituted $\stepsize(t)^{-1} = L + \extrastep(t)$ and then used the fact
that $\prox$ is strongly convex, so $\divergence(x(t + 1), x(t)) \ge \half
\norm{x(t) - x(t + 1)}^2$. By summing the
bound~\eqref{eqn:smooth-single-step-delayed-bound}, 
we have the following non-probabilistic inequality:
\begin{align}
  \lefteqn{\sum_{t=1}^T f(x(t + 1)) - f(x^*)} \nonumber \\
  & \le \<z(T+1), x(T + 1) - x^*\>
  + \frac{1}{\stepsize(T)} \prox(x(T + 1))
  + \sum_{t=1}^T \prox(x(t))\left[\frac{1}{\stepsize(t - 1)}
    - \frac{1}{\stepsize(t)}\right] \nonumber \\
  & ~~~~ - \sum_{t=1}^T \extrastep(t) \divergence(x(t + 1), x(t))
  + \sum_{t=1}^T \<\error(t), x(t + 1) - x^*\> \nonumber \\
  & \le \frac{1}{\stepsize(T + 1)} \prox(x^*)
  + \sum_{t=1}^T \prox(x(t)) \left[\frac{1}{\stepsize(t - 1)}
    - \frac{1}{\stepsize(t)}\right]
  - \sum_{t=1}^T \extrastep(t) \divergence(x(t + 1), x(t)) \nonumber \\
  & ~~~~
  + \sum_{t=1}^T \<\error(t), x(t + 1) - x^*\>
  \label{eqn:delayed-stoch-to-bound}
\end{align}
since $\prox(x) \ge 0$ and $x(T + 1)$ minimizes $\<z(T + 1), x\> +
\frac{1}{\stepsize(T + 1)} \prox(x)$.  What remains is to control the summed
$\error(t)$ terms in the bound~\eqref{eqn:delayed-stoch-to-bound}. We can do
this simply using the second part of
Lemma~\ref{lemma:bound-error-smooth}. Indeed, we have
\begin{align}
  \lefteqn{\sum_{t=1}^T \<\error(t), x(t + 1) - x^*\>}
  \label{eqn:split-stochastic-delayed-error} \\
  & = \sum_{t=1}^T \<\nabla f(x(t)) - \nabla
  f(x(t - \tau)), x(t + 1) - x^*\> + \sum_{t=1}^T
  \<\nabla f(x(t - \tau)) - g(t - \tau), x(t + 1) - x^*\>.
  \nonumber
\end{align}
We can apply Lemma~\ref{lemma:bound-error-smooth} to the first term in
\eqref{eqn:split-stochastic-delayed-error} by bounding $\norm{x(t) - x(t +
  1)}$ with Lemma~\ref{lemma:xt-difference-bound}. Since $\extrastep(t)
\propto \sqrt{t+\tau}$, Lemma~\ref{lemma:xt-difference-bound} with
$t_0 = \tau$ implies
$\E[\norm{x(t) - x(t + 1)}^2] \le \frac{4 G^2}{\extrastep(t)^2}$.  As a
consequence,
\begin{equation*}
  \E\bigg[\sum_{t=1}^T \<\nabla f(x(t)) - \nabla f(x(t - \tau)),
    x(t + 1) - x^*\>\bigg]
  \le 2 \tau G\radius
  + 2 L(\tau + 1)^2 G^2 \sum_{t=1}^T \frac{1}{\eta(t-\tau)^2}.
\end{equation*}

What remains, then, is to bound the stochastic (second) term in
\eqref{eqn:split-stochastic-delayed-error}. This is straightforward,
though:
\begin{align*}
  \lefteqn{\<\nabla f(x(t - \tau)) - g(t - \tau), x(t + 1) - x^*\>} \\
  & = \<\nabla f(x(t - \tau)) - g(t - \tau), x(t) - x^*\>
  + \<\nabla f(x(t - \tau)) - g(t - \tau), x(t + 1) - x(t)\> \\
  & \le \<\nabla f(x(t - \tau)) - g(t - \tau), x(t) - x^*\>
  + \frac{1}{2\extrastep(t)} \dnorm{\nabla f(x(t - \tau)) - g(t - \tau)}^2
  + \frac{\extrastep(t)}{2} \norm{x(t + 1) - x(t)}^2
\end{align*}
by the Fenchel-Young inequality applied to the conjugate pair $\half
\dnorm{\cdot}^2$ and $\half \norm{\cdot}^2$.  In addition, \mbox{$\nabla f(x(t
  - \tau)) - g(t - \tau)$} is independent of $x(t)$ given the sigma-field
containing $g(1), \ldots, g(t - \tau - 1)$, since $x(t)$ is a function of
gradients to time $t - \tau - 1$, so the first term has zero expectation. Also
recall that $\E[\dnorm{\nabla f(x(t-\tau)) - g(t-\tau)}]^2$ is
bounded by $\sigma^2$ by assumption. Combining the above two bounds into
\eqref{eqn:split-stochastic-delayed-error}, we see that
\begin{align}
  \lefteqn{\sum_{t=1}^T \E[\<\error(t), x(t + 1) - x^*\>]} \nonumber \\
  & \le \frac{\stddev^2}{2} \sum_{t=1}^T \frac{1}{\extrastep(t)}
  + \half \sum_{t=1}^T \extrastep(t) \norm{x(t + 1) - x(t)}^2
  + 2L G^2 (\tau + 1)^2
  \sum_{t=1}^T \frac{1}{\extrastep(t-\tau)^2} + 2 \tau G\radius.
  \label{eqn:clean-split-stochastic-delayed-error}
\end{align}

Since $\divergence(x(t + 1), x(t)) \ge \half \norm{x(t) - x(t + 1)}^2$,
combining \eqref{eqn:clean-split-stochastic-delayed-error} with
\eqref{eqn:delayed-stoch-to-bound} and noting that $\frac{1}{\stepsize(t -
  1)} - \frac{1}{\stepsize(t)} \le 0$ gives
\begin{equation*}
  \sum_{t=1}^T \E f(x(t + 1)) - f(x^*)
  \le \frac{1}{\stepsize(T + 1)} \prox(x^*)
  + \frac{\stddev^2}{2} \sum_{t=1}^T \frac{1}{\extrastep(t)}
  + 2 LG^2(\tau + 1)^2 \sum_{t=1}^T \frac{1}{\extrastep(t-\tau)^2}
  + 2 \tau G\radius.
\end{equation*}

\subsection{Proof of Theorem~\ref{theorem:md-delayed-stoch}}
\label{sec:md-stoch-delayed}

The proof of Theorem~\ref{theorem:md-delayed-stoch} is similar to that of
Theorem~\ref{theorem:da-delayed-stoch}, so we will be somewhat terse. We
define the error $\error(t) = \nabla f(x(t)) - g(t - \tau)$, identically as in
the earlier proof, and begin as we did in the proof of
Theorem~\ref{theorem:da-delayed-stoch}. Recall that
\begin{equation}
  f(x(t + 1)) - f(x^*)
  \le \<g(t - \tau), x(t + 1) - x^*\>
  + \<\error(t), x(t + 1) - x^*\> + \frac{L}{2} \norm{x(t) - x(t + 1)}^2.
  \label{eqn:first-order-and-lipschitz-ineq}
\end{equation}
Applying the first-order optimality condition to the definition of
$x(t+1)$~(\ref{eqn:md-update}), we get
\begin{equation*}
  \<\stepsize(t) g(t - \tau) + \nabla\prox(x(t + 1)) - \nabla \prox(x(t)),
  x - x(t + 1)\> \ge 0
\end{equation*}
for all $x \in \xdomain$. In particular, we have
\begin{align*}
  \stepsize(t)\<g(t - \tau), x(t + 1) - x^*\>
  & \le \<\nabla \prox(x(t + 1)) - \nabla \prox(x(t)), x^* - x(t + 1)\> \\
  & = \divergence(x^*, x(t)) - \divergence(x^*, x(t + 1)) -
  \divergence(x(t + 1), x(t)).
\end{align*}
Applying the above to the
inequality~\eqref{eqn:first-order-and-lipschitz-ineq}, we see that
\begin{align}
  \lefteqn{f(x(t + 1)) - f(x^*)} \nonumber \\
  & \le \frac{1}{\stepsize(t)}\left[
    \divergence(x^*, x(t)) - \divergence(x^*, x(t + 1))
    - \divergence(x(t + 1), x(t))\right] + \<\error(t), x(t + 1) -
  x^*\> + \frac{L}{2} \norm{x(t) - x(t + 1)}^2 
  \nonumber \\
  & \le \frac{1}{\stepsize(t)}\left[\divergence(x^*, x(t)) - \divergence(x^*,
    x(t + 1))\right]
  + \<\error(t), x(t + 1) - x^*\>
  - \extrastep(t) \divergence(x(t + 1), x(t))
  \label{eqn:md-single-step-bound}
\end{align}
where for the last inequality, we use the fact that $\divergence(x(t + 1),
x(t)) \ge \half \norm{x(t) - x(t + 1)}^2$, by the strong convexity of $\prox$,
and that $\stepsize(t)^{-1} = L + \extrastep(t)$.
By summing the inequality~\eqref{eqn:md-single-step-bound}, we have
\begin{align}
  \sum_{t=1}^T f(x(t + 1)) - f(x^*)
  & \le \frac{1}{\stepsize(1)} \divergence(x^*, x(1))
  + \sum_{t=2}^T \divergence(x^*, x(t)) \left[\frac{1}{\stepsize(t)}
    - \frac{1}{\stepsize(t - 1)}\right] \nonumber \\
  & \qquad ~ - \sum_{t=1}^T \extrastep(t) \divergence(x(t + 1), x(t)) 
  + \sum_{t=1}^T \<\error(t), x(t + 1) - x^*\>.
  \label{eqn:md-delayed-stoch-to-bound}
\end{align}
Comparing the bound~\eqref{eqn:md-delayed-stoch-to-bound} with the earlier
bound for the dual averaging algorithms~\eqref{eqn:delayed-stoch-to-bound}, we
see that the only essential difference is the $\stepsize(t)^{-1} -
\stepsize(t-1)^{-1}$ terms. The compactness assumption guarantees that
$\divergence(x^*, x(t)) \le \radius^2$, however, so
\begin{equation*}
  \sum_{t=2}^T \divergence(x^*, x(t))\left[\frac{1}{\stepsize(t)}
    - \frac{1}{\stepsize(t - 1)}\right]
  \le \frac{\radius^2}{\stepsize(T)}.
\end{equation*}
The remainder of the proof uses Lemmas~\ref{lemma:xt-difference-bound}
and~\ref{lemma:bound-error-smooth} completely identically to the proof of
Theorem~\ref{theorem:da-delayed-stoch}.

\subsection{Proof of Corollary~\ref{corollary:random-tau}}
\label{sec:corollary-simple-delayed}

We prove this result only for the mirror descent
algorithm~\eqref{eqn:md-delayed-update}, as the proof for the
dual-averaging-based algorithm~\eqref{eqn:da-delayed-update} is similar.  We
define the error at time $t$ to be $\error(t) = \nabla f(x(t)) - g(t -
\tau(t))$, and observe that we only need to control the second term
involving $\error(t)$ in the bound~\eqref{eqn:md-single-step-bound}
differently. Expanding the error terms above and using Fenchel's
inequality as in the proofs of Theorems~\ref{theorem:da-delayed-stoch}
and~\ref{theorem:md-delayed-stoch}, we have
\begin{align*}
  \lefteqn{\<\error(t), x(t + 1) - x^*\>} \\
  & \le \<\nabla f(x(t)) - \nabla f(x(t - \tau(t))), x(t + 1) - x^*\>
  + \<\nabla f(x(t - \tau(t))) - g(t - \tau(t)), x(t) - x^*\> \\
  & \qquad ~ + \frac{1}{2 \extrastep(t)} \dnorm{\nabla f(x(t - \tau(t)))
    - g(t - \tau(t))}^2 + \frac{\extrastep(t)}{2} \norm{x(t + 1) - x(t)}^2,
\end{align*}
Now we note that conditioned on the delay $\tau(t)$, we have
\begin{align*}
  \E[\norm{x(t - \tau(t)) - x(t + 1)}^2 \mid \tau(t)]
  & \le G^2 (\tau(t) + 1)^2 \stepsize(t - \tau(t))^2.
\end{align*}
Consequently we apply Lemma~\ref{lemma:bound-error-smooth} (specifically,
following the bounds~\eqref{eqn:bregmans-to-telescope}
and~\eqref{eqn:second-order-sum}) and find
\begin{align*}
  \lefteqn{\sum_{t=1}^T \<\nabla f(x(t)) - \nabla f(x(t - \tau(t))),
    x(t + 1) - x^*\>} \\
  & \le \sum_{t=1}^T
  \left[\divergencef(x^*, x(t)) - \divergencef(x^*, x(t - \tau(t)))\right]
  + G^2 \sum_{t=1}^T (\tau(t) + 1)^2 \stepsize(t - \tau(t))^2.
\end{align*}
The sum of $\divergencef$ terms telescopes, leaving only terms not received by
the gradient procedure within $T$ iterations, and we can use $\stepsize(t) \le
\frac{1}{\extrastep \sqrt{T}}$ for all $t$ to derive the further bound
\begin{equation}
  \label{eqn:how-many-divergences-appear}
  \sum_{t : t + \tau(t) > T} \divergencef(x^*, x(t))
  + \frac{G^2}{\extrastep^2 T} \sum_{t=1}^T (\tau(t) + 1)^2.
\end{equation}

To control the quantity~\eqref{eqn:how-many-divergences-appear}, all we need
is to bound the expected cardinality of the set $\{t \in [T] : t + \tau(t) >
T\}$.  Using Chebyshev's inequality and standard expectation bounds, we have
\begin{equation*}
  \E\left[\card(\{t \in [T] : t + \tau(t) > T\})\right]
  = \sum_{t=1}^T \P(t + \tau(t) > T)
  \le 1 + \sum_{t=1}^{T - 1} \frac{\E[\tau(t)^2]}{(T - t)^2}
  \le 1 + 2 B^2,
\end{equation*}
where the last inequality comes from our assumption that $\E[\tau(t)^2] \le
B^2$. As in Lemma~\ref{lemma:bound-error-smooth}, we have $\divergencef(x^*,
x(t)) \le 2G\radius$, which yields
\begin{equation*}
  \E\bigg[\sum_{t=1}^T\<\nabla f(x(t)) - \nabla f(x(t - \tau(t))),
    x(t + 1) - x^*\>\bigg]
  \le 6 G\radius B^2 + \frac{G^2 (B + 1)^2}{\extrastep^2}
\end{equation*}
We can control the remaining terms as in the proofs of
Theorems~\ref{theorem:da-delayed-stoch} and~\ref{theorem:md-delayed-stoch}.

\section{Proof of Theorem~\ref{theorem:da-distributed-delay}}
\label{sec:proof-da-distributed-delay}

The proof of Theorem~\ref{theorem:da-distributed-delay} is not too
difficult given our previous work---all we need to do is redefine the error
$\error(t)$ and use $\extrastep(t)$ to control the variance terms that arise.
To that end, we define the gradient error terms that we must
control. In this proof, we set
\begin{equation}
  \error(t) \defeq \nabla f(x(t)) - \sum_{i=1}^n \lambda_i g_i(t - \tau(i))
  \label{eqn:distributed-error}
\end{equation}
where $g_i(t) = \nabla f(x(t); \statsample_i(t))$ is the gradient of node $i$
computed at the parameter $x(t)$ and $\tau(i)$ is the delay
associated with node $i$. 

Using Assumption~\ref{assumption:lipschitz-grad} as in the proofs of
previous theorems, then applying Lemma~\ref{lemma:solution-convexity}, we
have
\begin{align*}
  & f(x(t + 1)) - f(x^*) \le \<\nabla f(x(t)), x(t + 1) - x^*\>
  + \frac{L}{2} \norm{x(t) - x(t + 1)}^2 \\
  & = \<\sum_{i=1}^n \lambda_i g_i(t - \tau(i)), x(t + 1) - x^*\>
  + \<\error(t), x(t + 1) - x^*\> + \frac{L}{2} \norm{x(t) - x(t + 1)}^2 \\
  & = \<z(t+1), x(t + 1) - x^*\> - \<z(t), x(t + 1) - x^*\>
  + \<\error(t), x(t + 1) - x^*\>
  + \frac{L}{2} \norm{x(t) - x(t + 1)}^2 \\
  & \le \<z(t+1), x(t + 1) - x^*\>
  - \<z(t), x(t) - x^*\> + \frac{1}{\stepsize(t)} \prox(x(t + 1))
  - \frac{1}{\stepsize(t)} \prox(x(t)) \\
  & ~~~~ - \frac{1}{\stepsize(t)}
  \divergence(x(t + 1), x(t)) + \<\error(t), x(t + 1) - x^*\>
  + \frac{L}{2} \norm{x(t) - x(t + 1)}^2.
\end{align*}
We telescope as in the proofs of Theorems~\ref{theorem:da-delayed-stoch}
and~\ref{theorem:md-delayed-stoch}, canceling $\frac{L}{2} \norm{x(t) - x(t +
  1)}^2$ with the $L \divergence$ divergence terms to see that
\begin{align}
  \lefteqn{\sum_{t=1}^T f(x(t + 1)) - f(x^*)} \nonumber \\
  & \le \<z(T+1), x(T + 1) - x^*\> + \frac{1}{\stepsize(T)} \prox(x(T))
  - \sum_{t=1}^T \extrastep (t) \divergence(x(t + 1), x(t))
  + \sum_{t=1}^T \<\error(t), x(t + 1) - x^*\> \nonumber \\
  & \le \frac{1}{\stepsize(T + 1)} \prox(x^*)
  - \sum_{t=1}^T \extrastep(t) \divergence(x(t + 1), x(t))
  + \sum_{t=1}^T \<\error(t), x(t + 1) - x^*\>.
  \label{eqn:distributed-stoch-to-bound}
\end{align}
This is exactly as in the non-probabilistic
bound~\eqref{eqn:delayed-stoch-to-bound} from the proof of
Theorem~\ref{theorem:da-delayed-stoch}, but the
definition~\eqref{eqn:distributed-error} of the error $\error(t)$ here is
different.

What remains is to control the error term in
\eqref{eqn:distributed-stoch-to-bound}. Writing the terms out, we have
\begin{align}
  \sum_{t=1}^T \<\error(t), x(t + 1) - x^*\>
  & = \sum_{t=1}^T \<\nabla f(x(t)) - \sum_{i=1}^n \lambda_i
  \nabla f(x(t - \tau(i))), x(t + 1) - x^*\> \nonumber \\
  & \qquad ~ + \sum_{t=1}^T \<\sum_{i=1}^n \lambda_i
  \left[\nabla f(x(t - \tau(i))) - g_i(t - \tau(i))\right],
  x(t + 1) - x^*\>
  \label{eqn:distributed-split-error}
\end{align}
Bounding the first term above
is simple via Lemma~\ref{lemma:bound-error-smooth}: as in the
proof of Theorem~\ref{theorem:da-delayed-stoch} earlier, we have
\begin{align*}
  \lefteqn{\E\bigg[\sum_{t=1}^T \<\nabla f(x(t)) - \sum_{i=1}^n \lambda_i
    \nabla f(x(t - \tau(i))), x(t + 1) - x^*\>\bigg]}\\
  & = \sum_{i=1}^n \lambda_i \sum_{t=1}^T
  \E[\<\nabla f(x(t)) - \nabla f(x(t - \tau(i))), x(t + 1) - x^*\>] \\
  & \le 2 \sum_{i=1}^n \lambda_i
  LG^2 (\tau(i) + 1)^2 \sum_{t=1}^T \frac{1}{\extrastep(t-\tau)^2}
  + \sum_{i=1}^n \lambda_i 2 \tau(i) G\radius.
\end{align*}

We use the same technique as the proof of
Theorem~\ref{theorem:da-delayed-stoch} to bound the second term from
\eqref{eqn:distributed-split-error}. Indeed, the
Fenchel-Young inequality gives
\begin{align*}
  \lefteqn{\<\sum_{i=1}^n \lambda_i
    \left[\nabla f(x(t - \tau(i))) - g_i(t - \tau(i))\right],
    x(t + 1) - x^*\>} \\
  & = \<\sum_{i=1}^n \lambda_i
  \left[\nabla f(x(t - \tau(i))) - g_i(t - \tau(i))\right],
  x(t) - x^*\> \\
  &\qquad +~ \<\sum_{i=1}^n \lambda_i
  \left[\nabla f(x(t - \tau(i))) - g_i(t - \tau(i))\right],
  x(t + 1) - x(t)\> \\
  & \le \<\sum_{i=1}^n \lambda_i
  \left[\nabla f(x(t - \tau(i))) - g_i(t - \tau(i))\right],
  x(t) - x^*\> \\
  & \qquad + ~ \frac{1}{2\extrastep(t)} \dnorm{
    \sum_{i=1}^n \lambda_i \left[\nabla f(x(t - \tau(i)))
      - g_i(t - \tau(i)) \right]}^2
  + \frac{\extrastep(t)}{2} \norm{x(t + 1) - x(t)}^2.
\end{align*}
By assumption, given the information at worker $i$ at time $t - \tau(i)$,
$g_i(t - \tau(i)))$ is independent of $x(t)$, so the first term has zero
expectation. More formally, this happens because $x(t)$ is a function of
gradients $g_i(1), \ldots, g_i(t - \tau(i) - 1)$ from each of the nodes $i$
and hence the expectation of the first term conditioned on $\{g_i(1), \ldots,
g_i(t - \tau(i) - 1)\}_{i=1}^n$ is 0.  The last term is canceled by the
Bregman divergence terms in \eqref{eqn:distributed-stoch-to-bound}, so
combining the bound~\eqref{eqn:distributed-split-error} with the above two
paragraphs yields
\begin{align*}
  \sum_{t=1}^T\E f(x(t + 1)) - f(x^*)
  & \le \frac{1}{\stepsize(t)} \prox(x^*)
  + 2 \sum_{i=1}^n \lambda_i LG^2 (\tau(i) + 1)^2
  \sum_{t=1}^T \frac{1}{\extrastep(t-\tau)^2}
  + \sum_{i=1}^n \lambda_i 2 \tau(i) G\radius \\
  & \qquad ~ + \sum_{t=1}^T \frac{1}{2\extrastep(t)}
  \E \bigg\|\sum_{i=1}^n \lambda_i \left[\nabla f(x(t - \tau(i)))
      - g_i(t - \tau(i))\right]\bigg\|_*^2.
\end{align*}

\section{Conclusion and Discussion}

In this paper, we have studied dual averaging and mirror descent algorithms
for smooth and non-smooth stochastic optimization in delayed settings, showing
applications of our results to distributed optimization.  We showed that for
smooth problems, we can preserve the performance benefits of parallelization
over centralized stochastic optimization even when we relax synchronization
requirements. Specifically, we presented methods that take advantage of
distributed computational resources and are robust to node failures,
communication latency, and node slowdowns. In addition, by distributing
computation for stochastic optimization problems, we were able to exploit
asynchronous processing without incurring any asymptotic penalty due to the
delays incurred. In addition, though we omit these results for brevity, it is
possible to extend all of our expected convergence results to guarantees with
high-probability.


\section*{Acknowledgments}

In performing this research, AA was supported by a Microsoft Research
Fellowship, and JCD was supported by the National Defense Science and
Engineering Graduate Fellowship (NDSEG) Program. We are very grateful to Ofer
Dekel, Ran Gilad-Bachrach, Ohad Shamir, and Lin Xiao for illuminating
conversations on distributed stochastic optimization and communication of
their proof of the bound~\eqref{eqn:smooth-da-stochastic}. We would also like
to thank Yoram Singer for reading a draft of this manuscript and giving useful
feedback.


\appendix

\section{Technical Results about Proximal Functions}
\label{sec:technical-prox}

In this section, we collect several useful results about proximal functions
and continuity properties of the solutions of proximal operators. We give
proofs of all uncited results in Appendix~\ref{appendix:prox-operators}.  We
begin with results useful for the dual-averaging updates~\eqref{eqn:da-update}
and~\eqref{eqn:da-delayed-update}.

We define the proximal dual function
\begin{equation}
  \label{eqn:prox-dual}
  \proxdual_\stepsize(z)
  \defeq \sup_{x \in \xdomain}
  \left\{\<-z, x\> - \frac{1}{\stepsize} \prox(x)\right\}.
\end{equation}
Since $\nabla \proxdual_\stepsize(z) = \argmax_{x \in \xdomain}
\{\<-z, x\> - \stepsize^{-1} \prox(x)\}$, it is clear
that $x(t) = \nabla \proxdual_{\stepsize(t)}(z(t))$. Further by strong
convexity of $\prox$, we have that
$\nabla\proxdual_\stepsize(z)$ is $\alpha$-Lipschitz continuous~\cite[Chapter
  X]{Nesterov09,HiriartUrrutyLe96b}, that is, for the norm $\norm{\cdot}$ with
respect to which $\prox$ is strongly convex and its associated dual norm
$\dnorm{\cdot}$,
\begin{equation}
  \label{eqn:lipschitz-dual}
  \norm{\nabla \proxdual_\stepsize(y) - \nabla \proxdual_\stepsize(z)}
  \le \stepsize \dnorm{y - z}.
\end{equation}
We will find one more result about solutions to the dual averaging update
useful. This result has essentially been proven in many
contexts~\cite{Nesterov09,Tseng08,DekelGiShXi10}.
\begin{lemma}
  \label{lemma:solution-convexity}
  Let $x^+$ minimize $\<z, x\> + A \prox(x)$ for all $x \in \xdomain$. Then
  for any $x \in \xdomain$,
  \begin{equation*}
    \<z, x\> + A \prox(x) \ge \<z, x^+\> + A \prox(x^+)
    + A \divergence(x, x^+)
  \end{equation*}
\end{lemma}

Now we turn to describing properties of the mirror-descent
step~\eqref{eqn:md-update}, which we will also use frequently.  The lemma
allows us to bound differences between $x(t)$ and $x(t + 1)$ for the
mirror-descent family of algorithms.
\begin{lemma}
  \label{lemma:md-closeness}
  Let $x^+$ minimize $\<g, x\> + \frac{1}{\stepsize} \divergence(x, y)$ over
  $x \in \xdomain$.  Then $\norm{x^+ - y} \le \stepsize\dnorm{g}$.
\end{lemma}

The last technical lemma we give explicitly bounds the differences between
$x(t)$ and $x(t + \tau)$, for some $\tau \ge 1$, by using the above continuity
lemmas.
\begin{lemma}
  \label{lemma:xt-difference-bound}
  Let Assumption~\ref{assumption:lipschitz-func} hold. Define $x(t)$ via the
  dual-averaging updates \eqref{eqn:da-update},~\eqref{eqn:da-delayed-update},
  or~\eqref{eqn:da-dist-update} or the mirror-descent
  updates~\eqref{eqn:md-update}, \eqref{eqn:md-delayed-update}, or
  \eqref{eqn:md-dist-update}. Let $\stepsize(t)^{-1} = L + \extrastep
  (t+t_0)^c$ for some $c \in [0, 1]$, $\extrastep > 0$, $t_0 \geq 0$, and $L
  \ge 0$. Then for any fixed $\tau$,
  \begin{equation*}
    \E[\norm{x(t) - x(t + \tau)}^2]
    \le \frac{4 G^2 \tau^2}{\extrastep^2 (t+t_0)^{2c}}
    ~~~ \mbox{and} ~~~
    \E[\norm{x(t) - x(t + \tau)}] \le \frac{2 G \tau}{\extrastep (t+t_0)^c}.
  \end{equation*}
\end{lemma}

\section{Proofs of Proximal Operator Properties}
\label{appendix:prox-operators}


\begin{proof-of-lemma}[\ref{lemma:md-closeness}]
  The inequality is clear when $x^+ = y$, so assume that $x^+ \neq y$.
  Since $x^+$ minimizes $\<g, x\> + \frac{1}{\stepsize} \divergence(x, y)$,
  the first order conditions for optimality imply
  \begin{equation*}
    \<\stepsize g + \nabla \prox(x^+) - \nabla\prox(y), x - x^+\> \ge 0
  \end{equation*}
  for any $x \in \xdomain$. Thus we can choose $y = x$ and see that
  \begin{equation*}
    \stepsize \<g, y - x\>
    \ge \<\nabla \prox(x^+) - \nabla \prox(y), x^+ - y\>
    \ge \norm{x^+ - y}^2,
  \end{equation*}
  where the last inequality follows from the strong convexity of $\prox$.
  Using H\"older's inequality gives that $\stepsize\dnorm{g}\norm{y - x} \ge
  \norm{x^+ - y}^2$, and dividing by $\norm{y - x}$ completes the proof.
\end{proof-of-lemma}

\begin{proof-of-lemma}[\ref{lemma:xt-difference-bound}]
  We first show the lemma for the dual-averaging updates.  Recall that $x(t) =
  \nabla \proxdual_{\stepsize(t)}(z(t))$ and $\nabla \proxdual_{\stepsize}$ is
  $\stepsize$-Lipschitz continuous. Using the triangle inequality,
  \begin{align}
    \norm{x(t) - x(t + \tau)}
    & = \norm{\nabla\proxdual_{\stepsize(t)}(z(t))
      - \nabla\proxdual_{\stepsize(t + \tau)}(z(t + \tau))}
    \nonumber \\
    & = \norm{\nabla\proxdual_{\stepsize(t)}(z(t)) -
      \nabla\proxdual_{\stepsize(t+\tau)}(z(t)) +
      \nabla\proxdual_{\stepsize(t+\tau)}(z(t)) -
      \nabla\proxdual_{\stepsize(t + \tau)}(z(t + \tau))}
    \nonumber \\
    & \le \norm{\nabla\proxdual_{\stepsize(t)}(z(t)) -
      \nabla\proxdual_{\stepsize(t+\tau)}(z(t))} +
    \norm{\nabla\proxdual_{\stepsize(t+\tau)}(z(t)) - 
      \nabla\proxdual_{\stepsize(t + \tau)}(z(t + \tau))}
    \nonumber \\
    & \le (\stepsize(t) - \stepsize(t + \tau)) \dnorm{z(t)} +
    \stepsize(t + \tau) \dnorm{z(t) - z(t + \tau)}. 
    \label{eqn:simple-difference-bound}
  \end{align}
  It is easy to check that for $c \in [0, 1]$,
  \begin{equation*}
    \stepsize(t) - \stepsize(t + \tau)
    \le \frac{c \extrastep \tau}{(L + \extrastep t^c)^2 t^{1-c}}
    \le \frac{c \tau}{\extrastep t^{1 + c}}.
  \end{equation*}

  By convexity of $\dnorm{\cdot}^2$, we can bound $\E[\dnorm{z(t) -
      z(t + \tau)}^2]$: 
  \begin{equation*}
    \E[\dnorm{z(t) - z(t + \tau)}^2]
    = \tau^2 \E\bigg[\bigg\|\frac{1}{\tau}\sum_{s=1}^\tau z(t + s)
        - z(t + s - 1)\bigg\|_*^2\bigg]
    = \tau^2 \E\bigg[\bigg\|\frac{1}{\tau}
      \sum_{s=0}^{\tau - 1} g(s)\bigg\|_*^2\bigg]
    \le \tau^2 G^2,
  \end{equation*}
  since $\E[\dnorm{\partial F(x; \statsample)}^2] \le G^2$ by assumption. Thus,
  bound~\eqref{eqn:simple-difference-bound} gives
  \begin{align*}
    \E[\norm{x(t) - x(t + \tau)}^2]
    & \le 2(\stepsize(t) - \stepsize(t + \tau))^2 \E[\dnorm{z(t)}^2]
    + 2 \stepsize(t + \tau)^2 \E[\dnorm{z(t) - z(t + \tau)}^2] \\
    & \le \frac{2 c^2 t^2 \tau^2 G^2}{\extrastep^2 t^{2 + 2c}}
    + 2 G^2 \tau^2 \stepsize(t + \tau)^2
    = \frac{2 c^2 \tau^2 G^2}{\extrastep^2 t^{2c}}
    + \frac{2 G^2 \tau^2}{(L + \extrastep(t + \tau)^c)^2},
  \end{align*}
  where we use Cauchy-Schwarz inequality in the first step. Since $c
  \le 1$, the last term is clearly bounded by $4 G^2 \tau^2 / 
  \extrastep^2 t^{2c}$.

  To get the slightly tighter bound on the first moment in the statement of
  the lemma, simply use the triangle inequality from the
  bound~\eqref{eqn:simple-difference-bound} and that $\sqrt{\E X^2}
  \ge \E |X|$. 

  The proof for the mirror-descent family of updates is similar. We focus on
  non-delayed update~\eqref{eqn:md-update}, as the other updates simply modify
  the indexing of $g(t + s)$ below. We know from
  Lemma~\ref{lemma:md-closeness} and the triangle inequality that
  \begin{equation*}
    \norm{x(t) - x(t + \tau)}
    \le \sum_{s = 1}^\tau \norm{x(t + s) - x(t + s - 1)}
    \le \sum_{s = 1}^\tau \stepsize(t + s - 1) \dnorm{g(t + s)}
  \end{equation*}
  Squaring the above bound, taking expectations, and recalling that
  $\stepsize(t)$ is non-increasing, we see
  \begin{align*}
    \E[\norm{x(t) - x(t + \tau)}^2]
    & \le \sum_{s = 1}^\tau \sum_{r = 1}^\tau \stepsize(t + s) \stepsize(t + r)
    \E[\dnorm{g(t + s)}\dnorm{g(t + r)}] \\
    & \le \tau^2\stepsize(t)^2
    \max_{r, s} \sqrt{\E[\dnorm{g(t + s)}^2]} \sqrt{\E[
        \dnorm{g(t + r)}^2]}
    \le \tau^2 \stepsize(t)^2 G^2
  \end{align*}
  by H\"older's inequality. Substituting the appropriate value for
  $\stepsize(t)$ completes the proof.
\end{proof-of-lemma}

\section{Error in~\cite{LangfordSmZi09}}
\label{appendix:langford-error}

Langford et al.~\cite{LangfordSmZi09}, in Lemma~1 of their paper, state an
upper bound on $\<g(t - \tau), x(t - \tau) - x^*\>$ that is essential to the
proofs of all of their results. However, the lemma only holds as an equality
for unconstrained optimization (i.e.\ when the set $\xdomain = \R^d$); in the
presence of constraints, it fails to hold (even as an upper bound). To see
why, we consider a simple one-dimensional example with $\xdomain = [-1,1]$,
$f(x) = |x|$, $\extrastep \equiv 2$ and we evaluate both sides of their lemma
with $\tau = 1$ and $t = 2$. The left hand side of their bound evaluates to
$1$, while the right hand side is $-1$, and the inequality claimed in the
lemma fails. The proofs of their main theorems rely on the application of
their Lemma~1 with equality, restricting those results only to unconstrained
optimization. However, the results also require boundedness of the gradients
$g(t)$ over all of $\xdomain$ as well as boundedness of the distance between
the iterates $x(t)$. Few convex functions have bounded gradients over all of
$\R^d$; and without constraints the iterates $x(t)$ are seldom bounded for all
iterations $t$.

\bibliographystyle{amsalpha}
\bibliography{bib}

\end{document}